%%%\magnification=\magstep1
\input amstex
\NoBlackBoxes
\define\II{\frak I}
\redefine\dd{\frak d}
\redefine\bb{\frak b}
\redefine\rr{{^{\omega}\omega}}
\define\stru{\langle H_{\aleph_2},{\in,}\omega_1, \frak I\rangle}
\define\rng{\text{\rm rng}}\define\xbody{\text{\rm body}}\define\Mi{\Bbb M}\define\lth{\text{\rm lth}}\define\nwd{\text{\rm NWD}}
\define\dom{\text{\rm dom}}\define\um{\Bbb U\Bbb M}

\documentstyle{amsppt}
\topmatter

\endtopmatter
\document
%\input prvni
%prvni
\
\vskip 2cm
\centerline{\bf CANONICAL MODELS FOR $\aleph_1$-COMBINATORICS\footnote{MR 
subject classification 03E50, 03E40. Secondary 03E35.}}
\vskip 1cm
\centerline{\smc Saharon Shelah\footnote{Research supported
by The Israel Science Foundation administered by The Israel
Academy of Sciences and Humanities. Publication number 610.}}
\smallskip
\centerline{Rutgers University}
\centerline{New Brunswick, NJ}
\smallskip

\centerline{Institute for Advanced Study}
\centerline{The Hebrew University}
\centerline{Jerusalem, Israel}
\vskip 1cm
\centerline{\smc Jind\v rich Zapletal\footnote{Research partially 
supported by NSF grant number DMS 9022140 and GA \v CR grant number 201/97/0216.}}
\smallskip
\centerline{Mathematical Sciences Research Institute}
\centerline{Berkeley, CA}
\vskip 2cm
\eightrm
{\bf Abstract.} We define the property of $\Pi_2$-compactness
of a statement $\phi$ of set theory, meaning roughly that the hard
core of the impact of $\phi$ on combinatorics of $\aleph_1$
can be isolated in a canonical model for the statement $\phi.$
We show that the following statements are $\Pi_2$-compact:
``dominating number$=\aleph_1,$''
``cofinality of the meager ideal$=\aleph_1$'', ``cofinality
of the null ideal$=\aleph_1$'',``bounding number$=\aleph_1$'', existence of various types
of Souslin trees and variations on uniformity of measure
and category$=\aleph_1.$  Several important new metamathematical
patterns among classical statements of set theory are
pointed out.
\tenrm 

\newpage
%\input I2
%i2
\head
{0. Introduction}
\endhead

One of the oldest enterprises in higher set theory is the study of 
combinatorics of the first uncountable cardinal. It appears that many phenomena
under investigation in this area are $\Sigma_2$ statements in
the structure $\langle H_{\aleph_2},\in,\II\rangle,$ where
$H_{\aleph_2}$ is the collection of sets of hereditary cardinality
$\aleph_1$ and $\II$ is a predicate for nonstationary subsets of $\omega_1.$
For example:

\roster
\item the Continuum Hypothesis--or, ``there exists an $\omega_1$ sequence of
reals such that every real appears on it'' \cite {Ca}
\item the negation of Souslin Hypothesis--or, ``there exists an $\omega_1$-tree
without an uncountable antichain'' \cite {So}
\item $\dd=\aleph_1$--or, ``there is a collection of $\aleph_1$ many functions
in $\rr$ such that any other such function is pointwise dominated
by one of them''
\item indeed, every equality $\frak x=\aleph_1$ for $\frak x$ a classical
invariant of the continuum is a $\Sigma_2$ statement--$\bb=\aleph_1,$
$\frak s=\aleph_1,$ additivity of measure$=\aleph_1\dots$ \cite {BJ}
\item there is a partition $h:[\omega_1]^2\to 2$ without an uncountable
homogeneous set \cite {T2}.
\item the nonstationary ideal is $\aleph_1$-dense \cite {W2}
\endroster

It appears that $\Sigma_2$ statements generally assert that the combinatorics
of $\aleph_1$ is complex. Therefore, given a sentence $\phi$ about sets,
it is interesting to look for models where $\phi$ and as few as possible
$\Sigma_2$ statements hold, in order to isolate the real impact of
$\phi$ to the combinatorics of $\aleph_1.$ The whole machinery of iterated
forcing \cite {S1} and more recently the $P_{max}$ method \cite {W2} 
were developed
explicitly for this purpose. This paper is devoted to constructing such
canonical $\Sigma_2$-poor (or $\Pi_2$-rich) models for a number
of classical statements $\phi.$

We consider cases of $\phi$ being $\dd=\aleph_1,$ cofinality of the
meager ideal$=\aleph_1,$ cofinality of the null ideal$=\aleph_1,$
$\bb=\aleph_1,$ existence of some variations of
Souslin trees, variations on uniformity of measure and category
$=\aleph_1$ and for all of these we find canonical models. It is also
proved that $\phi=$``reals can be covered by $\aleph_1$ many meager sets''
does not have such a model. But let us first
spell out exactly what makes our models canonical.

Fix a sentence $\phi.$ Following the $P_{max}$ method developed in \cite
{W2}, we shall aim for a $\sigma$-closed forcing $P_\phi$ definable in
$L(\Bbb R)$ so that the following holds:

\proclaim {Theorem Scheme 0.1}
Assume the Axiom of Determinacy in $L(\Bbb R).$ Then in $L(\Bbb R)^{P_\phi},$
the following holds:

\roster
\item ZFC, $\frak c=\aleph_2,$ the nonstationary ideal 
is saturated, $\delta_2^1=\aleph_2$
\item $\phi$
\endroster
\endproclaim

\proclaim {Theorem Scheme 0.2}
Assume that $\psi$ is a $\Pi_2$ statement for $\stru$ and

\roster
\item the Axiom of Determinacy holds in $L(\Bbb R)$
\item there is a Woodin cardinal with a measurable above it
\item $\phi$ holds
\item $\stru\models\psi.$
\endroster

Then in $L(\Bbb R)^{P_\phi},$ $\stru\models\psi.$
\endproclaim

If these two theorems can be proved for $\phi,$ we say that $\phi$
is $\Pi_2$-compact.

What exactly is going on? Recall that granted large cardinals, the theory
of $L(\Bbb R)$ is invariant under forcing \cite {W1}
and  so must be the theory of $L(\Bbb R)^{P_\phi}.$ Now varying
the ZFC universe enveloping $L(\Bbb R)$ so as to satisfy various
$\Pi_2$ statements $\psi,$ from Theorem Scheme 0.2 it follows that necessarily
$L(\Bbb R)^{P_\phi}$ must realize all such $\Pi_2$ sentences ever
achievable in conjunction with $\phi$ by forcing in presence of
large cardinals. In particular, roughly if $\psi_i:i\in I$ are
$\Pi_2$-sentences one by one consistent with $\phi$ then even their
conjunction is consistent with $\phi.$ And $L(\Bbb R)^{P_\phi}$
is {\it the} model isolating the impact of $\phi$ on combinatorics
of $\aleph_1.$ 

It is proved in \cite {W2} that $\phi=$``true'',``the nonstationary ideal 
is $\aleph_1$-dense'' and others are
$\Pi_2$-compact assertions. This paper provides many classical $\Sigma_2$
statements which are $\Pi_2$-compact as well as examples of natural
noncompact statements. In general, our results appear to run parallel with 
certain intuitions related to iterated forcing. The $\Pi_2$-compact
assertions often describe phenomena for which good preservation
theorems \cite {BJ, Chapter 6} are known. This is not surprising given
that in many cases the $P_{max}$ machinery can serve as a surrogate
to the preservation theorems--see Theorem 1.15(5)--and that many local
arguments in $P_{max}$ use classical forcing techniques--see Lemma 4.4
or Theorem 5.6. There are many open questions left:

\proclaim {Question 0.3}
Is it possible to define a similar notion of $\Pi_2$-compactness
without reference to large cardinals?
\endproclaim

\proclaim {Question 0.4}
Is the Continuum Hypothesis $\Pi_2$-compact? Of course, (1) of Theorem
Scheme 0.1 would have to be weakened to accomodate the Continuum Hypothesis.
\endproclaim

The first section outlines the proof scheme using which all the
$\Pi_2$ compactness results in this paper are demonstrated. The scheme
works subject to verification of three combinatorial properties
--Lemma schemes 1.10, 1.13, 1.16, of independent interest--of 
the statement $\phi$ in question, which is done in Sections 2-5.
These sections can be read and understood without any knowledge of
\cite {W2}. The only indispensable--and truly crucial-- reference
to \cite {W2} appears in the first line of the proof of Theorem 1.15. At the time this paper went into print, a draft version
of \cite {W2} could be obtained from its author. There
were cosmetical differences in the presentation of $P_{max}$ in this paper and in \cite {W2}. 

Our notation follows the set-theoretical standard as set forth in
\cite {J2}. The letter $\II$ stands for the nonstationary ideal
on $\omega_1.$
A system $a$ of countable sets is stationary if for every
function $f:(\bigcup a)^{<\omega}\to \bigcup a$ there is some $x\in
a$ closed under $f.$ $H_\kappa$ denotes the collection of all sets
of hereditary size $<\kappa.$ By a ``model'' we always mean a model
of ZFC if not explicitly said otherwise.
The symbol $\lozenge$ stands for the statement: there is
a sequence $\langle A_\alpha:\alpha\in\omega_1\rangle$ such that
$A_\alpha\subset\alpha$ for each $\alpha\in\omega_1$ and
for every $B\subset\omega_1$ the set $\{\alpha\in\omega_1:
B\cap\alpha=A_\alpha\}\subset\omega_1$ is stationary.
$\omega_1$-trees grow downward,
are always infinitely branching, are considered to consist
of functions from countable ordinals to $\omega$ ordered by reverse inclusion,
and if $T,\leq$ is such a tree then $T_\alpha=\{t\in T:$
ordertype of the set $\{ s\in T:t\leq s\}$ under $\geq$ is
just $\alpha\}$ and $T_{<\alpha}=\bigcup_{\beta\in\alpha} T_\beta.$
For $t\in T,$ $lev(t)$ is the unique ordinal $\alpha$ such that
$t\in T_\alpha.$ For trees $S$ of finite sequences, we write $[S]$
to mean the set $\{x:\forall n\in\omega\ x\restriction n\in S\}.$
When we compare open sets of reals sitting in different models
then we always mean to compare the open sets 
given by the respective definitions. $\delta_2^1$ is the supremum
of lengths of boldface $\Delta_2^1$ prewellorderings of reals,
 $\Theta$ is the supremum of lengths of all prewellorderings
of reals in $L(\Bbb R).$
In forcing, the western convention of writing $q\leq p$ if $q$
is more informative than $p$ is utilized. $\lessdot$ denotes
the relation of complete embedding between complete Boolean
algebras or partial orders. $RO(P)$ is the complete Boolean algebra 
determined by a partially ordered set $P,$
and $\Bbb C_\kappa$ is the Cohen algebra on $\kappa$ coordinates.
The algebra $\Bbb C=\Bbb C_\omega$ is construed as having a dense set
${^{<\omega}\omega}$ ordered by reverse inclusion.
Lemma and Theorem ``schemes'' indicate that we shall attempt
to prove some of their instances later.

The second author would like to thank Rutgers University for its 
hospitality in November 1995 and W. H. Woodin for many discussions
concerning his $P_{max}$ method.

%\input G2
%g2
\head
{1. General comments}
\endhead

This section sets up a framework in which all $\Pi_2$-compactness results 
in this paper will be proved. Subsection 1.0 introduces a crucial notion
of an iteration of a countable transitive model of ZFC. In Subsection 1.1,
a uniform in a sentence $\phi$ way of defining the forcing $P_\phi$
and proving instances of Theorem schemes 0.1, 0.2 is provided. In this 
{\it proof scheme} there are three combinatorial 
lemmas--1.10, 1.13, 1.16--which must
be demonstrated for each $\phi$ separately, and that is done in the section
of the paper dealing with that particular assertion $\phi.$ In Subsection
1.2 it is shown how subtle combinatorics of $\phi$ 
can yield regularity properties
of the forcing $P_\phi.$ And finally Subsection 1.3 gives some examples
of failure of $\Pi_2$-compactness.

\subhead
{1.0. Iterability}
\endsubhead

The cornerstone of the $P_{max}$ method is the possibility of finding
generic elementary embeddings of the universe with critical point equal
to $\omega_1.$ This can be done in several ways from sufficiently
large cardinals. Here is our choice:

\definition
{Definition 1.1} 
\cite {W1} Let $\delta$ be a Woodin cardinal. The nonstationary tower forcing
$\Bbb Q_{<\delta}$ is defined as the set 
$\{ a\in V_\delta:a$ is a stationary system
of countable sets$\}$ ordered by $b\leq a$ if for every $x\in b,$ $x\cap
\bigcup a\in a.$
\enddefinition

The important feature of this notion is the following. Whenever
$\delta$ is a Woodin cardinal and $G\subset\Bbb Q_{<\delta}$ is 
a generic filter then in $V[G]$ there is an ultrapower embedding
$j:V\to M$ such that the critical point of $j$ is $\omega_1^V,$
$j(\omega_1^V)=\delta$ and $M$ is closed under $\omega$ sequences;
in particular $M$ is wellfounded. All of this has been described and proved
in \cite {W1}. We shall be interested in iterations of this process.

\definition
{Definition 1.2}
\cite {W2}
Let $M$ be a countable transitive model of ZFC, $M\models
\delta$ is a Woodin cardinal. An iteration of
$M$  of length $\gamma$ based on $\delta$ is a sequence
$\langle M_\alpha:\alpha\in\gamma\rangle$
together with commuting maps $j_{\alpha\beta}:M_\alpha\to M_\beta:
\alpha\in\beta\in\gamma$ so that:

\roster
\item $M=M_0$
\item each $M_\alpha$ is a model of ZFC, possibly not transitive.
Moreover, $j_{\alpha\beta}$ are elementary embeddings.
\item for each $\alpha$ with $\alpha+1\in\gamma$ there is a 
$M_\alpha$ generic filter $G_\alpha\subset 
(\Bbb Q_{<j_{0\alpha}(\delta)})^{M_\alpha}.$
The model $M_{\alpha+1}$ is the generic ultrapower of $M_\alpha$     
by $G_\alpha$ and $j_{\alpha\alpha+1}$ is the ultrapower embedding.
\item at limit ordinals $\alpha\in\gamma$ a direct limit is taken.
\endroster
\enddefinition

\remark
{\bf Convention 1.3} If the models $M_\alpha$ are wellfounded 
we replace them with their transitive isomorphs.
Everywhere in this paper, in the context of one specific
iteration we keep the indexation system as in the above definition.
We write $\theta_\alpha=\omega_1^{M_\alpha}$ and $\Bbb Q_\alpha=
(\Bbb Q_{<j_{0\alpha}(\delta)})^{M_\alpha}.$
\endremark

\definition
{Definition 1.4} \cite {W2}
An iteration $j$ of a model $M$ is called full if
it is of length $\omega_1+1$ and for every pair $\langle x,\beta\rangle$
with $x\in \Bbb Q_\beta$ and $\beta\in\omega_1$ the set $\{\alpha\in
\omega_1:j_{\beta\alpha}(x)\in G_\alpha\}\subset\omega_1$ is
stationary.
\enddefinition

If all models in an iteration $j:M\to N$ of length $\omega_1+1$ are wellfounded
then $j$ can be thought of as stretching $\Cal P(\omega_1)^M$ into
a collection of subsets of the real $\omega_1.$ The fullness of $j$
is then a simple bookkeeping requirement on it, making sure in particular
that the model $N$ is correct about the nonstationary ideal, that is
$\frak I\cap N=\frak I^N.$

\definition
{Definition 1.5} \cite {W2}
A countable transitive model $M$ is said to be {\it iterable}
with respect to its Woodin cardinal $\delta$ if all of its iterations
based on $\delta$ produce only wellfounded models. $M$ is called
{\it stable iterable} with respect to $\delta$ if all of its generic
extensions by forcings of size $<\delta$ are iterable with
respect to $\delta.$
\enddefinition

It is not a priori clear whether iterability and stable iterability
are two different notions. We shall often neglect the dependence of the
above definition on the ordinal $\delta.$

Of course, a problem of great interest is to produce many rich iterable models.
The following lemma and its two corollaries record the two methods
of construction of such models used in this paper.

\proclaim {Lemma 1.6} \cite {W2}
Let $N$ be a transitive model of ZFC such that $\omega_1=On\cap N$ and
$N\models$``$\delta<\kappa$ are a Woodin and an inaccessible cardinal
respectively''. Then $M=N\cap V_\kappa$ is stable iterable
with respect to $\delta.$
\endproclaim

\proclaim {Corollary 1.7} \cite {W2}
Suppose that the Axiom of Determinacy holds in $L(\Bbb R).$ Then
for every real $x$ there is a stable iterable model containing $x.$
\endproclaim

\demo {Proof}
The determinacy assumption provides a model $N$ as in the Lemma
containing every real $x$ given beforehand \cite{Sc}.
 Then $x\in N\cap V_\kappa=M$
is the desired countable stable iterable model. \qed
\enddemo

\proclaim {Corollary 1.8} \cite {W2}
Suppose $\delta<\kappa$ are a Woodin and a measurable cardinal respectively.
Then for every real $x$ there is a stable iterable model elementarily
embeddable into $V_\kappa$ which contains $x.$
\endproclaim

\demo {Proof}
Fix a real $x$ and choose a countable elementary substructure $Z\prec
V_{\kappa+2}$ containing $x,\delta,\kappa$ and a measure $U$ on $\kappa.$
Let $\pi:Z\to\bar Z$ be the transitive collapse. Then the model
$\bar Z$ is iterable in Kunen's sense \cite {Ku} with respect to its measure
$\pi(U),$ since its iterations lift those of the universe using the measure
$U.$ Let $N^*$ be the $\omega_1$-th iterand of $\bar Z$
using the measure $\pi(U),$ let $N=N^*\cap V_{\omega_1}$ 
and let $M=N\cap V_{\pi(\kappa)}.$
The lemma applied to $N,\pi(\delta)$ and $\pi(\kappa)$ shows that
the model $M$ is stable iterable; moreover, $M=\bar Z\cap V_{\pi(\kappa)}$
and so the map $\pi\restriction M$ elementarily embeds $M$ into $V_\kappa.$
Since $x\in M,$ the proof is complete. \qed
\enddemo

\demo {Proof of Lemma}
We shall show that $M$ is iterable; the iterability of its small
generic extensions $M[G]$ follows from an application of our proof
to the model $N[G].$

For contradiction, assume that there is an iteration $j:M\to M^*$
which yields an illfounded model. Since $j''(On\cap M)$ is cofinal
in the ordinals of $M^*,$ there must be some $\beta\in M$ such that
$j(\beta)$ is illfounded. Choose the iteration $j$ of the minimal possible
length $\gamma_0$ and so that the least ordinal $\beta_0$ with
$j(\beta_0)$ illfounded is smallest possible among all iterations of length
$\gamma_0.$ Note that $\gamma_0$ must be a 
successor of a countable limit ordinal.

Now $\gamma_0,\beta_0$ are definable in the model $N$ as the unique
solutions to the formula $\psi(x,y,M)=$``for every large enough
cardinal $\lambda, Coll(\omega,\lambda)
\Vdash\chi(\check x,\check y,\check M)$'',
where $\chi(x,y,z)$ says ``$x$ is the minimal length of a bad iteration
of $z$ and $y$ is the minimal bad ordinal among such iterations of length
$x$''.  The point is that whenever $\kappa,\gamma,\beta<\lambda<\omega_1$
and $G\subset Coll(\omega,\lambda)$ is an $N$-generic filter, then
in the model $N[G]$ $\chi(\gamma_0,\beta_0,M)$ is a $\Sigma_2^1$ property
of hereditarily countable objects and therefore evaluated correctly.

There must be ordinals $\gamma_1<\gamma_0$ and
$\beta_1<j_{0\gamma_1}(\beta_0)$ such that $j_{\gamma_1\gamma_0}(\beta_1)$
is illfounded. Since $\kappa$ is an inaccessible cardinal of $N,$
the iteration $j_{0\gamma_1}$ can be copied to an iteration of the
whole model $N$ using the same nonstationary tower generic filters.
Write $j_{0,\gamma_1}:N\to N'$ for this extended version of $j_{0,\gamma_1}$
again and note that $M_{\gamma_1}=j_{0,\gamma_1}(M)=N'\cap V_{j_{0,\gamma_1}
(\kappa)}.$ By elementarity, the ordinals $j_{0,\gamma_1}(\gamma_0),
j_{0,\gamma_1}(\beta_0)$ are the unique solution to the formula
$\psi(x,y,j_{0,\gamma_1}(M))$ in the model $N'.$ However, an application
of the previous paragraph to $N'$ shows that this cannot be, since
$\gamma_0-\gamma_1,\beta_1$ are better candidates for such a solution.
Contradiction. \qed
\enddemo

\subhead
{1.1. The $P_{max}$ method}
\endsubhead

In this subsection we present a proof scheme used in this paper to show
that various $\Sigma_2$ sentences $\phi$ for the structure $\stru$
are $\Pi_2$ compact. For the record, all statements $\phi$ considered
here are consequences of $\lozenge$ and therefore easily found
consistent with large cardinals.

\definition
{Definition 1.9} 
The set $P_\phi$ is defined by induction on rank of its elements. $p\in P_\phi$
if $p=\langle M_p,w_p,\delta_p,H_p\rangle$ where if no confusion is
possible we drop the subscript $p$ and

\roster
\item $M$ is a countable transitive model of ZFC iterable with
respect to its Woodin cardinal $\delta$
\item $M\models w$ is a witness for $\phi$
\item $H\in M$ is the {\it history} of the condition $p;$ it is a set
(possibly empty) of pairs $\langle q,j\rangle$ where $q\in P_\phi$
and $j$ is in $M$ a full iteration of the model $M_q$ based on $\delta_q$
such that $j(w_q)=w$ and $j(H_q)\subset H$
\item if $\langle q,j\rangle,\langle q,k\rangle$ are both in $H$ then
$j=k.$
\endroster

The ordering on $P_\phi$ is defined by $q\leq p$ just in case $\langle
p,j\rangle\in H_q$ for some $j.$
\enddefinition

The notion of a witness for $\phi$ used above is the natural one; if
$\phi=\exists x\ \forall y\ \chi(x,y)$ with $\chi$ a $\Sigma_0$ formula,
then $x\in H_{\aleph_2}$ is a witness for $\phi$ whenever 
$\stru\models\forall y\ \chi(x,y).$ However, for obviously equivalent
versions of the sentence $\phi$ this notion can vary a little.
A special care will always be taken as to what variation of $\phi$
we are working with.

The idea behind the definition of the forcing $P_\phi$ is to construct
$H_{\aleph_2}$ of the resulting model as a sort of direct limit
of its approximations in countable models taken under iterations
--which are recorded in the histories--and extensions. 

The possibility of $\Pi_2$-compactness of $\phi$ depends on the validity
of three combinatorial lemmas which show how witnesses for $\phi$
in countable transitive models can be stretched by iterations
of these models into real witnesses for $\phi.$ These lemmas are used
in Theorem 1.15 for $\sigma$-closure and various 
density arguments about $P_\phi.$

The first combinatorial fact to be proved is:

\proclaim {Lemma scheme 1.10} (Simple Iteration Lemma)
Suppose $\lozenge$ holds. If $M$ is a countable transitive model
of ZFC iterable with respect to its Woodin cardinal $\delta$
and $M\models$``$w$ is a witness for $\phi$'' then there is a
full iteration $j$ based on $\delta$
of the model $M$ such that $j(w)$ is a witness for $\phi.$
\endproclaim

Certainly there is a need for some assumption of the order of $\lozenge,$
since a priori $\phi$ does not have to hold at all and then $j(w)$ could not be a witness for it! Later we shall try
to optimalise this assumption to $\stru\models\phi,$ the weakest
possible.

For a detailed analysis of the forcing $P_\phi$ a more involved variant
of this lemma will be necessary. Essentially, the iteration $j$ is to be built
cooperatively by two players, one of whom attempts to make $j(w)$
into a witness for $\phi.$ The other one stages various
local obstacles to that goal. The relevant definitions:

\definition {Definition 1.11}
A sequence $\vec N$ of models with a witness 
is a system $\langle w,N_i,\delta_i:
i\in\omega\rangle$ where

\roster
\item $N_i$ are countable transitive models of ZFC+$\delta_i$ is a Woodin
cardinal+$w$ is a witness for $\phi;$ we set $\Bbb Q_i=\Bbb Q_{<\delta_i}$
as computed in $N_i$
\item $N_i\in N_{i+1}$ and $\omega_1^{N_i}$ is the same for all $i\in\omega$
\item if $N_i\models$``$a\in V_{\delta_i}$ is a stationary system of
coutable sets'' then $N_{i+1}\models$``$a$ is a stationary system of countable
sets''; so $\Bbb Q_0\subset\Bbb Q_1\subset\dots$
\item if $N_i\models$``$A\subset\Bbb Q_i$ is a maximal antichain''
then $N_{i+1}\models$``$A\subset\Bbb Q_{i+1}$ 
is a maximal antichain''.
\endroster

We say that the sequence begins with the triple $\langle N_0,w,\delta_0\rangle,$
set $\Bbb Q_{\vec N}=\bigcup_{i\in\omega}\Bbb Q_i$
and $\omega_1^{\vec N}=\omega_1^{N_0}.$ A filter $G\subset
\Bbb Q_{\vec N}$ is said to be $\vec N$ generic if it meets
all maximal antichains of $\Bbb Q_{\vec N}$ which happen  to belong
to $\bigcup_{i\in\omega}N_i.$
\enddefinition

This definition may seem a little artificial, an artifact of the machinery of
\cite {W2}. The really interesting information a sequence of models
carries is the model $\bigcup_{i\in\omega}N_i$ with its first-order theory.
This model can be viewed as a $w$-correct extension of $N_0.$
It is important that

\roster
\item $\langle\bigcup_i N_i,\in\rangle\models w$ is a witness for $\phi$
\item  $\langle\bigcup_i N_i,\in\rangle$ satisfies all $\Pi_2$-consequences of
ZFC in the language $\in,\frak S,$ where $\frak S$ is the predicate for
stationary systems of countable sets
\item $N_0$ is in $\bigcup_i N_i$ correct about stationary systems
of countable sets in $V_{\delta_0}$ and their maximal antichains.
\endroster

It should be noted that though $\Bbb Q_{\vec N}$ is not an element
of $\bigcup N_i,$ it is a class in that model--the class of all
stationary systems of countable sets. If a filter $G\subset\Bbb Q_{\vec N}$
is $\vec N$-generic then the filters $G\cap\Bbb Q_i$ are $N_i$-generic
by (3,4) of the above definition. However, not every $N_i$-generic
filter on $\Bbb Q_i$ can be extended into an $\vec N$-generic
filter on $\Bbb Q_{\vec N}.$

\definition {Definition 1.12}
$\Cal G_\phi$ is a two-person game of length $\omega_1$ between players
Good and Bad. The rules are:

{\rm Round 0:} The player Bad plays $M,w,\delta$ such that $M$ is
a countable transitive model of ZFC iterable with respect to its
Woodin cardinal $\delta$ and $M\models$``$w$ is a witness for $\phi$''

{\rm Round $\alpha>0$:} an ordinal $\gamma_\alpha$ and an iteration
$j_{\gamma_\alpha}:M\to M_{\gamma_\alpha}$ of length $\gamma_\alpha+1$
based on $\delta$ are given.
\roster
\item"{$\bullet$}"Bad plays a sequence $\vec N$ of models beginning with
$M_{\gamma_\alpha},j_{\gamma_\alpha}(w),
j_{\gamma_\alpha}(\delta)$ and a condition $p\in\Bbb Q_{\vec N}.$
\item"{$\bullet$}"Good plays an 
$\vec N$-generic filter $G\subset \Bbb Q_{\vec N}$
with $p\in G.$
\item"{$\bullet$}"Bad plays an ordinal $\gamma_{\alpha+1}>\gamma_\alpha$
and an iteration $j_{\gamma_{\alpha+1}}$ of $M$ of length $\gamma_{\alpha+1}
+1$ which prolongs the iteration $j_{\gamma_\alpha}$ and such that
the $\gamma_\alpha$-th ultrapower on it is taken using the filter
$G\cap\Bbb Q_{\gamma_\alpha}.$
\endroster

Here, $\gamma_1=-1, j_1=id$ and at limit $\alpha$'s, $j_{\gamma_\alpha}$ is the
direct limit of the iterations played before $\alpha.$

In the end, let $j$ be the direct limit of the iterations played. The
player Good wins if either the player Bad cannot play at some stage
or the iteration $j$ is not full or $j(w)$ is a witness for $\phi.$
\enddefinition

Thus the player Bad is responsible for the bookkeeping to make the iteration
full and has a great freedom in prolonging the iteration
on a nonstationary set of steps. The player Good has a limited access
on a closed unbounded set of steps to steering $j(w)$ into a witness for
$\phi.$ In the real life, the player Bad can easily play all the way
through $\omega_1$ and make the resulting iteration full.

We shall want to prove

\proclaim {Lemma scheme 1.13}
(Strategic Iteration Lemma) Suppose $\lozenge$ holds. Then the player Good has 
a winning strategy in the game $\Cal G_\phi.$
\endproclaim

Now suppose that the relevant instances of Lemma schemes 1.10, 1.13 are true for
$\phi.$ Then, granted the Axiom of Determinacy in $L(\Bbb R),$
the model $L(\Bbb R)^{P_\phi}$ can be completely analysed using the
methods of \cite {W2} to verify Theorem scheme 0.1.
Let $G\subset P_\phi$ be a generic filter.

\definition {Definition 1.14}
In $L(\Bbb R)[G],$ for any $p\in G$ define

\roster
\item $k_p$ is the iteration of $M_p$ which is the direct limit of the system
$\{ j:\exists q\in G\ \langle p,j\rangle\in H_q\}.$
\item $W=k_p(w_p).$
\endroster
\enddefinition

It is obvious from the definition of the poset $P_\phi$ that the
system $\{ j:\exists q\in G\ \langle p,j\rangle\in H_q\}$ is directed,
and that the definition of $W$ does not depend on the particular choice
of $p\in G.$

\proclaim {Theorem 1.15}
Assume the Axiom of Determinacy in $L(\Bbb R)$ and the relevant
instances of Lemma schemes 1.10 and 1.13 hold. Then 
$P_\phi$ is a $\sigma$-closed notion of forcing and in $L(\Bbb R)[G],$
the following hold:

\roster
\item ZFC
\item for every $X\in H_{\aleph_2}$ there are $p\in G$ and $x\in M_p$
such that $X=k_p(x).$
\item $\frak c=\aleph_2,$ the nonstationary ideal 
is saturated, $\delta_2^1=\omega_2$
\item $\stru\models\phi$ and $W$ is a witness for $\phi$
\item Suppose $\psi=\forall x\ \exists y\ \chi(x,y)$ for some
$\Sigma_0$-formula $\chi$ is a $\Pi_2$-statement for $\stru.$
Suppose that ZFC+$\lozenge$ proves that for each $x\in H_{\aleph_2}$
there is a forcing $P$ preserving witnesses to $\phi$ and stationary subsets
of $\omega_1$ such that $P\Vdash\exists y\chi(\check x,\dot y).$ Then
$\stru\models\psi.$
\endroster
\endproclaim

\demo {Proof}
Parts (1,2,3) are straigthforward generalizations of Section 4.3 in \cite {W2}.
Work in $L(\Bbb R)[G]$ and prove (4).

First note that whenever $p\in G$ then $k_p''(\II)^{M_p}=\II\cap
k_p''M_p.$ To see it, suppose $p\in G,M_p\models$``$s\subset\omega_1$
is a stationary set'' and fix a club $C\subset\omega_1.$ By (2),
there is a condition $q\in G$ and $c\in M_q$ so that $C=k_q(c).$
Let $r\in G$ be a common lower bound of $p,q$ with $\langle p,i\rangle,
\langle q,j\rangle\in H_r.$ Then

\roster
\item $M_r\models j(c)\subset\omega_1$ is a club
\item $M_r\models i(s)$ is stationary, since the iteration $i$ is full in $M_r.$
\endroster

Therefore $j(c)\cap i(s)\neq 0.$ By absoluteness, $k_rj(c)\cap k_ri(s)\neq 0$
and since $k_rj=k_q$ and $k_ri=k_p,$ we have $C\cap k_p(s)=0$ and $k_p(s)$
is stationary.

Now suppose (4) fails; so $\phi=\exists x\ \forall y\ \upsilon(x,y)$
for some $\Sigma_0$ formula $\upsilon,$ and $W$ is not a witness for
$\phi,$ and $\lnot\upsilon(W,Y)$ for some $Y\in H_{\aleph_2}.$ By (2),
there is a condition $p\in G$ and $y\in M_p$ such that $Y=k_p(y).$

But now, $M_p\models\stru\models\upsilon(w_p,y),$ since $w_p$
is a witness for $\phi$ in the model $M_p.$ By elementarity
of $k_p,$ absoluteness and the previous paragraph $\stru\models
\upsilon(j_p(w_p)=W,j_p(y)=Y),$ contradiction.

To prove (5), fix the formula $\chi$ and note that by (2) it is enough to
show that for $p\in P_\phi$ and $x\in M_p$ there is $q\leq p$
which forces an existence of $Y$ such that $\stru\models\chi
(k_p(x),Y).$ And indeed, using Corollary 1.7 choose a countable transitive 
stable iterable model
$M$ with $M\models$``$\lozenge+\delta$ is a Woodin cardinal''
with $p\in M.$ Apply the Iteration Lemma 1.10 
in $M$ to get a full iteration
$j$ of $M_p$ such that $j(w_p)$ is a witness for $\phi$ in $M.$
By the assumptions on $\chi$ applied in $M\cap V_\kappa,$ where $\kappa$
is the least inaccessible cardinal of $M,$ there is a generic extension
$M[K]$ of $M$ by a forcing of size $<\kappa$ 
preserving $j(w_p)$ and stationary subsets 
of $\omega_1$ such that there is $y\in M[K]$ with
$\chi(j(x),y).$ Obviously, setting $q=\langle M[K],j(w_p),\delta, H
\rangle,$ where $H=j(H_p)\cup\{\langle p,j\rangle\},$
we have $q\leq p$ and $q\Vdash\chi(\dot k_p(x),k_q(y))$ as desired. \qed
\enddemo

The rudimentary comparison of the cardinal structure of $L(\Bbb R)$
and $L(\Bbb R)[G]$ carries over literally from \cite {W2}: namely
$\aleph_1,\aleph_2$ are the same in these models, $\Theta=\aleph_3^{
L(\Bbb R)[G]}$ and all cardinals above $\Theta$ are preserved.
This will not be used anywhere in this paper.

It should be remarked that under the assumptions of the Theorem,
the model $K=L(\Cal P(\omega_1))$ as evaluated in $L(\Bbb R)[G]$ satisfies 
(1)-(5) and it can be argued that $K$ is the ``canonical model''
in view of its minimal form and Theorem 1.23. In fact, $L(\Bbb R)[G]$
is a generic extension of $K$ by the poset $(\omega_2^{<\omega_2})^K.$

The final point in the analysis of the model $L(\Bbb R)[G]$ is the proof
of Theorem scheme 0.2 for $P_\phi.$ We know of only one approach for doing this,
namely to prove

\proclaim {Lemma scheme 1.16}
(Optimal Iteration Lemma) Suppose $\stru\models\phi.$ 
Whenever $M$ is a countable transitive model iterable with respect
to its Woodin cardinal $\delta$ and $M\models$``
$w$ is a witness for $\phi$'' then there is a full iteration $j$ of $M$
based on $\delta$ such that $j(w)$ is a witness for $\phi.$
\endproclaim
  
It is crucial that the assumption of this optimal iteration lemma
is truly the weakest possible. Provided Lemma schemes 1.10, 1.13
and 1.16 are true for
$\phi,$ we can conclude

\proclaim {Corollary 1.17}
Suppose instances of Iteration Lemmas 1.10, 1.13 and 1.16 for $\phi$
are true. Then $\phi$ is $\Pi_2$-compact.
\endproclaim

\demo {Proof}
We shall prove the relevant instance of Theorem Scheme 0.2. Assume that $\psi$
is a $\Pi_2$ sentence, $\psi=\forall x\ \exists y\ \chi(x,y)$
for some $\Sigma_0$ formula $\chi.$ Assume that there is a Woodin
cardinal $\delta$ with a measurable cardinal $\kappa$ above it, and
$\stru\models\psi\land\phi.$ It must be proved that $L(\Bbb R)^{P_\phi}
\models\stru\models\psi.$

For contradiction suppose that $p\in P_\phi$ forces $\lnot\psi=\exists
x\ \forall y\ \lnot \chi(x,y)$ holds in $\stru$ of the generic extension.
By Theorem 1.15 (2), by eventually strengthening the condition $p$
we may assume that there is $x\in\stru^{M_p}$ so that $p\Vdash
\forall y\lnot\chi(k_p(x),y)$ holds in $\stru$ of the generic extension.

Following Corollary 1.8, there is a countable transitive iterable model
$M$ elementarily embeddable into $V_\kappa$ containing $p,$
which is a hereditarily countable object in $M.$
By the relevant instance of the Optimal Iteration Lemma applied within $M$
there is a full iteration $j$ of the model $M_p$
such that $j(w_p)$ is a witness for $\phi$ in $\stru^M.$
It follows that the quadruple $q=\langle M,j(w_p),\bar \delta,H\rangle,$
where $H=j(H_p)\cup\{\langle p,j\rangle\}$ and 
$\bar\delta$ is a Woodin cardinal of $M,$ is a condition in $P_\phi$
and $q\leq p.$ Since in $M,$ $\stru\models\psi,$ necessarily there is
$y\in \stru^M$ such that $\stru^M\models\chi(j(x),y).$ It follows that
$q\Vdash\chi(k_qj(x)=k_p(x),k_q(y))$ in $\stru$ of the generic extension,
a contradiction to our assumptions on $p,x.$ \qed
\enddemo

Predicates other than $\frak I$ can be added to the language
of $H_{\aleph_2}$ keeping the amended version of Theorem scheme 0.2 valid.
For example, if Iteration Lemmas 1.10, 1.13 and 1.16 hold for $\phi$
then the relevant instance of Theorem scheme 0.2 can be shown to hold
with the richer structure $\langle H_{\aleph_2},\in,\frak I,X:X\subset\Bbb R,
X\in L(\Bbb R)\rangle;$ however, the proof is a little involved
and we omit it. See \cite {W2}. In certain cases a predicate for witnesses
for $\phi$ can be added keeping Theorem scheme 0.2
true for $\phi.$ This increases the expressive power of the language
a little. Such a possibility will be discussed on a case-by-case basis.

Let us recapitulate what we proved in this subsection. 
Let $\phi$ be a $\Sigma_2$ sentence for $\stru,$ a consequence of
$\lozenge.$ If instances of
Lemma schemes 1.10 and 1.13 
are shown to hold, then the model $L(\Bbb R)^{P_\phi}$
has the properties listed in Theorem 1.15
or Theorem scheme 0.1. And if the optimal iteration lemma 1.16
for $\phi$ is proved then the relevant instance of Theorem scheme 0.2
is true and $\phi$ is $\Pi_2$-compact. It should be noted that
Iteration Lemma 1.10 follows from both Lemma 1.13 and Lemma 1.16.
We include it because
it is frequently much easier to prove and because it is often the first
indication that a $\Pi_2$-compactness type of result can be proved.

\subhead
{1.2. Order of witnesses}
\endsubhead

After an inspection of the proofs of iteration lemmas in the subsequent
sections the following notion comes to light:

\definition
{Definition 1.18} Let $\phi$ be a $\Sigma_2$ sentence. For $v,w\in H_{\aleph_2}$
we set $v\leq_\phi w$ if in every forcing extension of the universe
whenever $v$ is a witness to $\phi$ then $w$ is such a witness.
\enddefinition

Of course, the formally impermissible consideration of all forcing extensions
can be expressed as quantification over partially ordered sets. While
restricting ourselves to just {\it forcing} extensions may seem to be
somewhat artificial, it is logically the easiest way and the resulting
notion fits all the needs of the present paper. It should be noted
that $\leq_\phi$ is sensitive to the exact definition of a witness as it
was the case for $P_\phi.$ Obviously $\leq_\phi$ is a quasiorder
and the nonwitnesses form the $\leq_\phi$ smallest
$\leq_\phi$-equivalence class.

\remark {\bf Example 1.19}
For $\phi=$``there is a Souslin tree'' and $S,T$ such trees the relation
$T\leq_\phi S$ is equivalent to the assertion ``for every $s\in S$ there are
$s'\leq s,t\in T$ such that $RO(S\restriction s')$ can be completely
embedded into $RO(T\restriction t)$''. For then, preservation of the Souslinity
of $T$ implies the preservation of the Souslinity of $S.$ On the other
hand, if the assertion fails, there must be $s\in S$ such that
$T\Vdash$``$S\restriction s$ is an Aronszajn tree'' because every cofinal branch
through $S$ is generic. By the c.c.c. productivity theorem then,
the finite condition forcing specializing the tree $S\restriction s$
preserves the Souslinity of $T$ and collapses the Souslinity of $S;$
ergo, $T\not\leq_\phi S.$
\endremark

Note that in the above example the relation $\leq_\phi$ was $\Sigma_1$
on the set of all witnesses. It is not clear whether this behavior
is typical; the proofs of the iteration lemmas in this paper
always use a $\Sigma_1$ phenomenon to guarantee the relation $\leq_\phi$
or the $\leq_\phi$-equivalence of two witnesses.

\definition
{Definition 1.20} Suppose $\phi=\exists x\ \forall y\ \chi(x,y)$ for some
$\Sigma_0$ formula $\chi$ and let $\psi(x_0,x_1)$ be $\Sigma_1.$
We say that $\psi$ is a {\it copying procedure} for $\phi$ if
ZFC $\vdash\stru\models\forall x_0,x_1\ (\psi(x_0,x_1)\to(\forall y\ 
\chi(x_0,y)\leftrightarrow\forall y\ \chi(x_1,y))).$ In other words,
$\psi(x_0,x_1)$ guarantees that $x_1$ is a witness for $\phi$
iff $x_0$ is a witness for $\phi.$
\enddefinition

\remark
{\bf Example 1.21} Let $\phi=$``there is a nonmeager set of reals of size
$\aleph_1$''. One possible copying procedure 
for $\phi$ is $\psi(x_0,x_1)=$``there
is a continuous category-preserving function $f:\Bbb R\to\Bbb R$ such that
$f''(x_0)=x_1$''. Note that this is really a statement about a code
for $f,$ which is essentially a real, and it can be cast in a $\Sigma_1$ form.
\endremark

The following theorems, quoted without proof,
 are applications of the above concepts.
The first implies that the forcings $P_\phi$ for sentences $\phi$ considered 
in this paper are
all homogeneous and therefore the $\Sigma_n$ theory of $L(\Bbb R)^{P_\phi}$ is a 
(definable) element of $L(\Bbb R)$ for every $n\in\omega.$
The second shows that models in sections 2,3 and 5 not only optimalize
the $\Sigma_2$-theory of $\stru$ but are in fact characterized
by this property.
 The choice of copying procedures necessary
for its proof will always be clear from the arguments in the section
dealing with that particular $\phi.$

\proclaim {Theorem 1.22}
Suppose the axiom of determinacy holds in $L(\Bbb R),$ suppose $\phi$
is a $\Sigma_2$ sentence for which iteration lemmas 1.10, 1.13 hold and $\psi$
is a copying procedure for $\phi$ such that ZFC proves
one of the following:

\roster
\item for every two witnesses $x_0,x_1$ for $\phi$ there is a forcing $P$
preserving stationary subsets of $\omega_1$ and witnesses for $\phi$
such that $P\Vdash\psi(x_0,x_1)$
\item for every witness $x_0$ and every countable transitive iterable
model $M$ with $M\models$``$w$ is a witness for $\phi$'' there
is a full iteration $j$ of $M$ such that $\psi(x_0,j(w))$ holds.
\endroster

Then $P_\phi$ is a homogeneous notion of forcing.
\endproclaim

\proclaim {Theorem 1.23}
Suppose the axiom of determinacy holds in $L(\Bbb R),$
suppose $\phi$ is a $\Sigma_2$ sentence for which 
iteration lemmas 1.10, 1.13 hold and $\psi$ is a copying procedure
for $\phi$ such that ZFC proves ``for every two witnesses $x_0,x_1$
to $\phi$ there is a forcing $P$ preserving stationary subsets of $\omega_1$
and witnesses to $\phi$ such that $P\Vdash\psi(x_0,x_1).$
{\bf If} the $\Sigma_2$-theory of the structure
$\langle H_{\aleph_2},\in,\frak I,X:X\subset\Bbb R,X\in L(\Bbb R)\rangle$
is the same in $V$ as in $L(\Bbb R)^{P_\phi}$ {\bf then} $\Cal P(\omega_1)
=\Cal P(\omega_1)\cap L(\Bbb R)[G]$ for some possibly external
$L(\Bbb R)$-generic filter $G\subset P_\phi.$ 
\endproclaim

\subhead
{1.3. Limitations}
\endsubhead

Of course by far not every $\Sigma_2$ sentence $\phi$ can be handled using
the proof scheme outlined in Subsection 1.1. Each of the three iteration lemmas
can prove to be a problem; in some cases, 
it is possible to show that the statement
$\phi$ is not $\Pi_2$-compact by exhibiting $\Pi_2$ assertions $\psi_i:
i\in I$ each of whom is consistent with $\phi$ yet $\bigwedge_{i\in I}
\psi_i\vdash\lnot\phi.$

\remark
{\bf Example 1.24} 
The simple iteration lemma for $\phi=$``the Continuum Hypothesis''
fails. The reason is that whenever $M$ is a countable transitive
iterable model and $j$ is an iteration of $M$ then $j(\Bbb R\cap M)\neq\Bbb R$
--namely, the real coding the model $M$ is missing from $j(\Bbb R
\cap M).$
\endremark

\remark
{\bf Example 1.25} The simple iteration lemma for $\phi=$``there is a maximal almost
disjoint family (MAD) of sets of integers of size $\aleph_1$'' fails.
Note that if $A\subset\Cal P({^{<\omega}\omega})$ is a MAD extending the
set of all branches in ${^{<\omega}\omega}$ then $A$ is collapsed as a MAD
whenever a new real is added to the universe. Thus if $M$ is a countable
transitive iterable model with $M\models$``the Continuum Hypothesis holds
and $A$ is a MAD as above'' then no iteration $j$ of $M$ makes $j(A)$
into a MAD for the same reason in the previous example.
\endremark

\remark
{\bf Example 1.26}
 The strategic iteration lemma for $\phi=$``there is a nonmeager
set of reals of size $\aleph_1$'' with the natural notion of witness 
fails. The reason is somewhat arcane and we omit it.
\endremark

\remark
{\bf Example 1.27} The optimal iteration lemma for $\phi=$``the reals can be
covered with $\aleph_1$ many meager sets'' cannot be proved. For
suppose that $M$ is a countable transitive iterable model and 
$M\models$``the Continuum Hypothesis holds and $\Cal C=\{ X_f:f\in\rr\},$
where $X_f\subset\rr$ is the set of all reals pointwise dominated
by the function $f.$ So $\Cal C$ constitutes a covering of the real
line by $\aleph_1$ meager sets.'' Also suppose that
$\phi\land\bb>\aleph_1$ holds in the universe--this is consistent
and happens after adding $\aleph_2$ Laver reals to a model of GCH \cite {Lv}.
Then no iteration $j$ of the model $M$ can make $j(\Cal C)$ into a covering
of the real line, because there always will be a function in $\rr$
eventually dominating all of $j(\Bbb R\cap M).$ Note that in this
case, $\Cal C$ should be thought of as a collection of Borel codes
as opposed to a set of sets of reals.
\endremark

The previous example suggests that $\phi$ is not $\Pi_2$-compact,
and indeed, it is not. For consider $\Pi_2$ sentences

$$\eqalign{\psi_0=&\text{``}\bb>\aleph_1",\cr
\psi_1=&\text{``for every bounded family }A\subset\rr\text{ of size }
\aleph_1\cr
\ &\text{ there is a function infinitely many times 
equal to every function in }A
".\cr}$$

Now $\phi\land\psi_0$ holds after iterating Laver reals, 
$\phi\land\psi_1$ holds after iterating proper $\rr$-bounding
forcings \cite {S3, Proposition 2.10}, 
and $\psi_0\land\psi_1\vdash\lnot\phi$ can be
derived easily from the combinatorial characterization of $\phi$
in \cite {Ba, BJ}..

\remark
{\bf Example 1.28} The optimal iteration lemma for $\phi=$``$\frak t=\aleph_1$''
cannot be proved. Recall that $\frak t$ is the minimal length of a tower
and a tower is a decreasing sequence of infinite subsets of $\omega$
without lower bound in the modulo finite inclusion ordering. To see
the reason for the failure, suppose $M$ is a countable transitive iterable
model and $M\models$``the Continuum Hypothesis holds and $t$ is a tower
of height $\omega_1$ consisting of sets of asymptotic density one.''
That such towers exist under CH has been pointed out to us by W. Hugh Woodin.
Suppose that in the universe $\frak t=\aleph_1$ holds and no towers consist
of sets from a fixed Borel filter--such a situation can be attained
by iterating Souslin c.c.c. forcings over a model of CH. Then no iteration
$j$ of $M$ can make $j(t)$ into a tower.
\endremark

It seems that the two $\Pi_2$ assertions

$$\eqalign
{\psi_0&=\text{``no towers consist of sets from a fixed Borel filter''}\cr
\psi_1&=\text{``every tower consists of sets from some Borel filter''}\cr}$$

provide a witness for non-$\Pi_2$-compactness of $\phi,$
however, the consistency of $\phi\land\psi_1$ seems to be a difficult open
problem.

\remark
{\bf Example 1.29} The optimal iteration lemma for $\phi=$``there
is a Souslin tree'' cannot be proved. For suppose that
$M$ is a countable transitive iterable model with $M\models$
``$\lozenge$ and $T$ is a homogeneous Souslin tree''. Suppose
that in the universe there are Souslin trees but none of them
are homogeneous--this was proved consistent in \cite {AS}.
Then no iteration $j$ of $M$ can make $j(T)$ into a Souslin tree,
since $j(T)$ is necessarily homogeneous.
\endremark

Again, the above example provides natural candidates to witness the
non-$\Pi_2$-compactness of $\phi.$ Let

$$\eqalign{\psi_0=&\text{``for every finite set }
T_i:i\in I\text{ of Souslin trees 
there are }t_i\in T_i\cr
\ &\text{ such that }\prod_iT_i\restriction t_i
\text{ is c.c.c.''},\cr
\psi_1=&\text{``for every Souslin tree $T$ there are finitely many }
t_i:i\in I\text{ in }T\cr
\ &\text{ such that }\prod_iT\restriction t_i\text
{ is nowhere c.c.c.''}.\cr}$$

The sentence $\phi\land\psi_0$ was found consistent by \cite {AS}, but the consistency of $\phi\land\psi_1$ is an open problem. In view of
the results of Subsection 4.0 the sentences $\psi_0,\psi_1$ are the {\it only}
candidates for noncompactness of $\phi.$

%\input D2
%d2
\head
{2. Dominating number}
\endhead

The proof of $\Pi_2$-compactness of the sentence ``there is a family
of $\aleph_1$ many functions in $\rr$ such that any function in $\rr$
is modulo finite dominated by one in the family'' or $\dd=\aleph_1,$
is in some sense prototypical, and the argument will be adapted to
other invariants in Section 3. The important concept we isolate to
prove the iteration lemmas is that of {\it subgenericity}. It essentially
states that the classical Hechler forcing is the optimal way for adding
a dominating real. To our knowledge, this concept has not been explicitly
defined before.

\subhead
{2.0. The combinatorics of $\dd$}
\endsubhead

 There is a natural Souslin \cite
{JS} forcing associated to the order of eventual dominance on $\rr$
designed to add a ``large'' function:

\definition {Definition 2.1}
The Hechler forcing $\Bbb D$ is the set 
$\{\langle a,A\rangle:\dom(a)=n$ for some $n\in\omega,$ $\rng(a)\subset
\omega$ and $A$ is a finite subset of $\rr\}.$ The order is defined
by $\langle a,A\rangle\leq\langle b,B\rangle$ if
\roster
\item $b\subset a, B\subset A$
\item $\forall n\in\dom(a\setminus b)\ \forall f\in B\ a(n)\geq f(n).$
\endroster

For a condition $p=\langle a,A\rangle$ the function $\xbody(p)\in\rr$
is defined as $\xbody(p)(n)=a(n)$ if $n\in\dom(a)$ and $\xbody(p)(n)=
\max\{ f(n):f\in A\}$ if $n\notin dom(a).$

If $G\subset \Bbb D$ is a generic filter, the Hechler real $d$
is defined as $\bigcup\{a:\langle a,0\rangle\in G\}.$
\enddefinition

Below, we shall make use of
restricted versions of $\Bbb D.$ Say $f\in\rr.$ Then define $\Bbb D(f)$ to
be the set of all $p\in \Bbb D$ with $\xbody(p)$ pointwise dominated
by the function $f,$ with the order inherited from $\Bbb D.$ Note that
$\Bbb D$ as defined above is not a separative poset.

Obviously, all forcings defined above are c.c.c. The important combinatorial
fact about  $\Bbb D$ is that a Hechler real is in fact an ``optimal''
function eventually dominating every ground model function.
 This will be immediately made precise:

\definition {Definition 2.2}
Let $M$ be a transitive model of ZFC and $f\in\rr.$ We say that
the function $f$ $\Bbb D$-dominates $M$ if every $g\in M\cap\rr$
is eventually dominated by $f.$
\enddefinition

\proclaim {Lemma 2.3}
Let $M$ be a transitive model of ZFC and let $f\in\rr$ $\Bbb D$-dominate
$M.$ If $D\subset \Bbb D\cap M$ is a dense set in $M,$ then
$D\cap \Bbb D(f)$ is dense in $\Bbb D(f)\cap M.$
\endproclaim

Note that $\Bbb D\cap M$ is $\Bbb D$ as computed in $M;$ 
also $\Bbb D(f)\cap M\notin M.$

\demo {Proof}
Fix a dense set $D\subset \Bbb D\cap M$ in $M$ and a condition 
$p=\langle a,A\rangle$ in $\Bbb D(f)\cap M.$ We shall produce a condition 
$q\leq p$ in $D\cap \Bbb D(f).$ Working in $M,$ it is easy to construct
a sequence $\langle a_i,A_i\rangle:i\in\omega$ of conditions in $D$
with $\dom(a_i)=n_i$ and

\roster
\item $\langle a_i,A_i\rangle\leq p$
\item $\langle a_{i+1},A_{i+1}\rangle$ is any element of $D$ stronger than
$\langle\xbody(p)\restriction n_i,A_i\rangle.$
\endroster

Define a function $g\in\rr$ as follows: for $n\leq n_0,$ let $g(n)=a(n).$
For $n\geq n_0,$ find an integer $i\in\omega$ with $n_i\leq n<n_{i+1}$
and set $g(n)=\max\{ a_{i+1}(n),h(n):h\in A_i\}.$ Since $g\in\rr\cap M,$
the function $f$ dominates $g$ pointwise starting from some $n_i.$
Then $q=\langle a_{i+1},A_{i+1}\rangle\in D\cap \Bbb D(f)$ is the desired
condition. \qed
\enddemo

\proclaim {Corollary 2.4}
(Subgenericity) Let $P$ be a forcing and $\dot g$ a $P$-name such that
\roster
\item $P\Vdash\dot g\in\rr$ $\Bbb D$-dominates the ground model
\item for every $f\in\rr$ the boolean value $||\check f\leq \dot g$ pointwise
$||_P$ is non-zero.
\endroster
Then there is a complete embedding $RO(\Bbb D)\lessdot P*(\Bbb D(\dot g)\cap$
the ground model$)=R$ so that $R\Vdash$``$\dot d\leq\dot g$ pointwise'',
where $\dot d$ is the $\Bbb D$-generic real.
\endproclaim

Thus under every $\Bbb D$-dominating real, a Hechler real
is lurking behind the scenes.

\demo {Proof}
The Hechler real $\dot d$ will be read off the second
iterand in the natural way, and Lemma 2.3 will guarantee its genericity. 
By some  Boolean algebra theory, that yields a complete embedding
of $RO(\Bbb D\restriction p)$ to $RO(R),$ for some $0\neq p\in\Bbb D.$
We must prove that $p=1.$ But fix an arbitrary $q=\langle a,A\rangle
\in \Bbb D$ and let $g=\xbody(q).$ Then $\langle ||\check g\leq\dot f$
pointwise$||_P,\check q\rangle\in R$ is a nonzero element
of $RO(R)$ forcing that the $\Bbb D$-generic real meets 
the condition $q.$ Thus with the given embedding, any condition in
$\Bbb D$ can be met, consequently $p=1$ and the proof is
complete. \qed
\enddemo

\proclaim {Corolary 2.5}
Let $M$ be a countable transitive model of ZFC such that $M\models$
``$P$ is a forcing adding a dominating function $\dot f$
as in Corollary 2.4'', let $p\in P$ and 
let $g\in\rr$ be any function, not necessarily
in the model $M.$ Then there is an $M$-generic filter $G\subset B$
containing $p$ such that $g$ is eventually dominated by $\dot f/G.$
\endproclaim

\demo {Proof}
Apply the previous Corollary in the model $M$ and find the forcing $R$
and the relevant embeddings $P\lessdot R,$ $\Bbb D\lessdot R.$ 
Assume $P\subset R$ and set $x$ to be the projection of $r$ into $\Bbb D$
via the above embedding. Step out of the model $M$
 and find a filter $H\subset\Bbb D$ such that

\roster
\item $x\in H$
\item $H$ meets every maximal antichain of $\Bbb D$ which is an element of
$M$
\item the Hechler real $e$ 
derived from $H$ eventually dominates the function $g.$
\endroster

This is easily done. Now the 
key point is that the model $M$ computes maximal antichains of $\Bbb D$
correctly: if $M\models$``$A\subset \Bbb D$ is a maximal antichain''
then this is a $\Pi^1_1$ fact about $A$ under suitable coding
and therefore $A$ really is a maximal antichain of $\Bbb D.$ Consequently,
the filter $H\cap M\subset\Bbb D^M$ is $M$-generic.

Choose an $M$-generic filter $K
\subset R$ with $H\cap M\subset K$ under the embedding of $\Bbb D$
mentioned above. Let $G=K\cap P.$ The filter $G\subset P$
is $M$-generic and has the desired properties: the function $g$ is eventually
dominated by $e$ which is pointwise smaller than $\dot f/G.$  \qed
\enddemo

\subhead
{2.1. A model for $\frak d=\aleph_1$}
\endsubhead

A natural notion of a witness for $\dd=\aleph_1$ to be used in the definition
of $P_{\dd=\aleph_1}$ is that of an eventual domination cofinal subset of $\rr$
of size $\aleph_1.$ We like to consider an innocent strengthening of this
notion in order to later ensure that assumptions of Corollary 2.5 are satisfied.

\proclaim {Lemma 2.6}
The following are equivalent.
\roster 
\item $\frak d=\aleph_1$
\item there is a sequence $d:\omega_1\to\rr$ increasing in 
the eventual domination order 
such that for every
$f\in\rr$ the set $S_f=\{\alpha\in\omega_1:f\leq d(\alpha)$
pointwise$\}$ is stationary.
\endroster
\endproclaim

A sequence $d$ as in (2) will be called a {\it good} dominating
sequence and will be used as a witness for $\dd=\aleph_1.$

\demo {Proof}
Only (1)$\to$(2) needs an argument. Choose an arbitrary 
eventual domination cofinal
set $\{ f_\alpha:\alpha\in\omega_1\}\subset\rr$ and a sequence
$\langle S_{a,\alpha}:a\in{^{<\omega}\omega},\alpha\in\omega_1\rangle$
of pairwise disjoint stationary subsets of $\omega_1.$ By a
straightforward induction on $\beta\in\omega_1$ it is easy to
build the sequence $d:\omega_1\to\rr$ so that

\roster
\item $d(\beta)$ eventually dominates all $f_\alpha:\alpha\in\beta$
and $d(\alpha):\alpha\in\beta$
\item if $\beta\in S_{a,\alpha}$ then $d(\beta)$ pointwise dominates
both $a$ and $f_\alpha.$
\endroster

The sequence $d$ is as required. For choose $f\in\rr.$ There are
$a\in{^{<\omega}\omega},\alpha\in\omega_1$ so that $f$ is pointwise
dominated by the function $g$ taking maxima of functional values
of $a$ and $f_\alpha.$ The set $S_f$ is then a superset of $S_{a,\alpha}$
and therefore stationary.\qed
\enddemo

Suppose now that $M$ is a countable transitive model of ZFC,
$M\models$``$d:\omega_1\to\rr$ is a good dominating sequence
and $\delta$ is a Woodin cardinal''. 
 Working in $M$ a simple observation is that
$\Bbb Q_{<\delta}\Vdash$``$(j_{\Bbb Q}d)(\omega_1^M)\ 
\Bbb D\text{-dominates }M$'',
where $j_{\Bbb Q}$ is the term for the generic nonstationary tower embedding.
Moreover, by (2) above the pair $\Bbb Q_{<\delta}, (j_{\Bbb Q}d)(\omega_1^M)$
satisfies requirements of Corollary 2.5. Thus there is a generic
ultrapower of $M$ lifting $(j_{\Bbb Q}d)(\omega_1^M)$ arbitrarily
high in the eventual domination 
order in $V.$ Also, whenever $j$ is a full iteration
of the model $M$ such that $j(d)$ is a dominating sequence,
it is really a good dominating sequence.

\proclaim {Optimal Iteration Lemma 2.7}
 Suppose $\dd=\aleph_1.$ Whenever $M$ is a countable transitive
model iterable with respect to its Woodin cardinal $\delta$ and
$M\models$``$d$ is a good dominating sequence'' there is a full
iteration $j$ of $M$ so that $j(d)$ is a good dominating sequence.
\endproclaim

\demo {Proof}
Let $\{ f_\alpha:\alpha\in\omega_1\}$ be an eventual domination 
cofinal family of
functions. We shall produce a full iteration
$j$ of the model $M$ based on $\delta$ 
with $\theta_\alpha=\omega_1^{M_\alpha}$ such
that the function $jd(\theta_\alpha)$ eventually dominates the function
$f_\alpha,$ for every $\alpha\in\omega_1.$ This will prove the lemma.

The iteration $j$ will be constructed by induction on $\alpha\in\omega_1.$
First, fix a partition $\{ S_\xi:\xi\in\omega_1\}$ of
$\omega_1$ into pairwise disjoint stationary sets. By induction
on $\alpha\in\omega_1$ build models $M_\alpha$ together with the elementary
embeddings plus an enumeration $\{ \langle x_\xi,\beta_\xi\rangle:
\xi\in\omega_1\}$ of all pairs $\langle x,\beta\rangle$ with
$x\in \Bbb Q_\beta.$ The induction hypotheses at $\alpha$ are:

\roster
\item the function $jd(\theta_\alpha)$ eventually dominates $f_\alpha$
\item the initial segment $\{ \langle x_\xi,\beta_\xi\rangle:
\xi\in\theta_\alpha\}$ of the enumeration under construction has 
been built and it enumerates all pairs $\langle x,\beta\rangle$ with
$x\in \Bbb Q_\beta,$ $\beta\in\alpha.$
\item for $\gamma\in\alpha$ if  $\theta_\gamma\in S_\xi$ for some
(unique) $\xi\in\theta_\gamma$ then $j_{\beta_\xi,\gamma}(x_\xi)\in
G_\gamma.$
\endroster

The hypothesis (1) ensures that the resulting sequence $jd$ will
be dominating. The enumeration together with (3) will imply
the fullness of the iteration.

At limit stages, the direct limit of the previous models and the union 
of the enumerations constructed so far is taken. The successor step is
handled easily using a version of Corollary 2.5 
below the condition $j_{\beta_\xi,\gamma}(x_\xi)\in \Bbb Q_\alpha$
if $\alpha\in S_\xi$ for some $\xi\in\theta_\alpha$
and the observation just before
the formulation of Lemma 2.7. Let $j$ be the direct limit
of the iteration system constructed in the induction process. \qed
\enddemo

This lemma could have been proved even without subgenericity,
since after all, the forcing $\Bbb Q$ adds a Hechler real by design.
With the sequences of models entering the stage though it is important
to have a sort of a uniform term for this real.

\proclaim {Strategic Iteration Lemma 2.8} Suppose $\dd=\aleph_1.$ 
The player Good has a winning
strategy in the game $\Cal G_{\dd=\aleph_1}.$
\endproclaim

\demo {Proof}
Let $\vec N=\langle d,N_i,\delta_i:i\in\omega\rangle$ be a sequence
of models with a good dominating sequence $d,$ let $y_0\in\Bbb Q_{
\vec N}$ and let $f\in\rr$ be an arbitrary function. We shall show that
there is an $\vec N$-generic filter $G\subset\Bbb Q_{\vec N}$ with $y_0\in G$
such that $(j_\Bbb Qd)(\omega_1^{\vec N})$ eventually dominates
the function $f,$ where $j_\Bbb Q$ is the generic ultrapower 
embedding of the model $N_0$ using the filter $G\cap \Bbb Q_{<\delta_0}^{N_0}.$
With this fact a winning strategy in the game $\Cal G_{\dd=\aleph_1}$
for the good player consists of an appropriate bookkeeping using
a fixed dominating sequence as in the previous lemma.

First, let us fix some notation. Choose an integer $i\in\omega,$ work in
$N_i$ and set $\Bbb Q_i=\Bbb Q_{<\delta_i}^{N_i}.$ Consider the
$\Bbb Q_i$-term $j_i$ for a $\Bbb Q_i$-generic ultrapower embedding
of the model $N_i.$ The function $j_id(\omega_1^{\vec N})$ is forced
to be represented by the function $\alpha\mapsto d(\alpha)$ and to
$\Bbb D$-dominate the model $N_i.$ Applying Corollary 2.4 in $N_i$
to $\Bbb Q_i$ and $j_id(\omega_1^{\vec N})$ it is possible to 
choose a particular dense subset $R_i$ of the iteration found in that Corollary,
namely $R_i=\{\langle y,a,A\rangle:y\in\Bbb Q_i,\langle a,A\rangle\in\Bbb D,
\ A\subset\rng(d)$ and for every $x\in y$ the function $d(x\cap\omega_1)$
pointwise dominates $\xbody(\langle a,A\rangle)\}$ ordered by
$\langle z,b,B\rangle\leq\langle y,a,A\rangle$ just in case
$z\leq y$ in $\Bbb Q_i$ and $\langle b,B\rangle\leq\langle a,A\rangle$
in $\Bbb D.$ It is possible to restrict ourselves to sets
$A\subset\rng(d)$ since $rng(d)$ is an eventual domination cofinal family
in $N_i.$ As in that corollary, the $\Bbb Q_i$-generic will be read
off the first coordinate and the $\Bbb D$-generic real $e\in\rr$
will be read off the other two, with $j_id(\omega_1)$ pointwise
dominating the function $e.$ With this embedding of $\Bbb D$ into
the poset $R_i,$ we can compute the projection
$pr_{i\Bbb D}(\langle y,a,A\rangle)=\Sigma_{\Bbb D}\{\langle b,B\rangle\in
\Bbb D:\langle b,B\rangle\leq\langle a,A\rangle,B\subset\rng(d)$ and
the system $z=\{x\in y:d(x\cap\omega_1)$ pointwise dominates
the function $\xbody(\langle b,B\rangle)\}$ is stationary$\}.$

Now step out of the model $N_i.$
There are two key points, capturing the uniformity of the above definitions
in $i\in\omega:$

\roster
\item $R_0\subset R_1\subset\dots$
\item $pr_{i\Bbb D}(\langle y,a,A\rangle)$ computes the same value
in $\Bbb D$ in all models $N_i$ with $y\in\Bbb Q_i.$
\endroster

Therefore we can write $pr_{\Bbb D}(\langle y,a,A\rangle)$ to mean
the constant value in $\Bbb D$ of this expression without any
danger of confusion. Another formulation of (2) is that $N_0$
computes a function from $\Bbb D$ into $R_0$ which constitutes a complete
embedding of $\Bbb D$ into all $R_i$ in the respective models $N_i.$
Note that $R_0$ is {\it not} a complete suborder of the $R_i$'s.

Now everything is ready to construct the filter $G\subset \Bbb Q_{\vec N}.$
First, let us choose a sufficiently generic filter $H\subset\Bbb D.$
There are the following requirements on $H:$

\roster
\item"{(3)}" $pr_{\Bbb D}(\langle y_0,0,0\rangle)\in H$
\item"{(4)}" $H$ meets all maximal antichains of $\Bbb D$ which happen to be 
elements of $\bigcup_i N_i$
\item"{(5)}" the Hechler real $e\in\rr$ given by the filter $H$ 
eventually dominates the function $f\in\rr.$
\endroster

This is easily arranged. It follows from (4) that $H\cap N_i$ is an 
$N_i$-generic subset of $\Bbb D^{N_i}$ since the model $N_i$ is $\Sigma_1^1$
correct and therefore computes maximal antichains of $\Bbb D$ correctly.

Let $X_k:k\in\omega$ be an enumeration of all maximal antichains of
$\Bbb Q_{\vec N}$ which are elements of $\bigcup_i N_i.$ By induction
on $k\in\omega$ build a descending sequence $y_0\geq y_1\geq\dots
\geq y_k\geq\dots$ of conditions in $\Bbb Q_{\vec N}$ so that

\roster
\item"{(6)}" $y_{k+1}$ has an element of $X_k$ above it
\item"{(7)}" $pr_{\Bbb D}(\langle y_{k+1},0,0\rangle)\in H$
\endroster

This is possible by the genericity of the filter $H.$ Suppose $y_k$ is given.
There is an integer $i\in\omega$ such that $y_k\in \Bbb Q_i$ and
$X_k\in N_i$ is a maximal antichain in $\Bbb Q_i.$ Now $H\cap N_i$
is a Hechler $N_i$-generic filter and $\bar X_k=\{
\langle z,0,0\rangle:z\in X_k\}\subset R_i$ is a maximal antichain
in $R_i.$ Therefore, there must be a condition $y_{k+1}\leq y_k$
as required in (6,7).

In the end, let $G\subset\Bbb Q_{\vec N}$ be the filter generated 
by the conditions $y_k:k\in\omega.$ It is an $\vec N$-generic filter
 by (6) and the function $j_0d(\omega_1^{\vec N})$ pointwise dominates
the Hechler function $e$ by (7) 
and therefore--from (5)--eventually dominates the
function $f\in\rr$ as desired. \qed
\enddemo

\proclaim {Conclusion 2.11}
The sentence $\phi=\dd=\aleph_1$ is $\Pi_2$-compact, moreover, in
Theorem Scheme 0.2 we can add a predicate for dominating sequences
of length $\omega_1$ to the language of $\stru.$
\endproclaim

\demo {Proof}
All the necessary iteration lemmas have been proved. To see that
the dominating predicate $\frak D$ can be added, go through the proof
of Corollary 1.15 again and note that if $j$ is a full iteration of
a countable transitive iterable model $M$ such that $j(M\cap\Bbb R)$
is cofinal in the eventual domination ordering then $\frak D\cap M_{\omega_1}=
\frak D^{M_{\omega_1}}.$ \qed
\enddemo 

%\input O2
%o2
\head
{3. Other $\dd$-like cardinal invariants}
\endhead

The behavior of the dominating number seems to be typical for
a number of other cardinal invariants. We present here two cases
which can be analysed completely. Recall that for an arbitrary
ideal, the cofinality of that ideal is defined as the minimal size
of a collection of small sets such that any small set is covered
by one in that collection. This is an important cardinal characteristic
of that ideal \cite {BJ}.

\subhead
{3.0. Cofinality of the meager ideal}
\endsubhead

In this subsection we prove that the statement ``the cofinality of the meager
ideal is $\aleph_1$'' is $\Pi_2$-compact. It is not difficult to see and 
will be proved below that it is enough to pay attention to the
nowhere dense ideal. As in the previous section, there is a canonical
forcing related to this ideal.

\definition {Definition 3.1}
\roster
\item $\nwd$ is the set of all perfect nowhere dense trees on $\{0,1\}$
\item for a tree $T\in\nwd$ and 
a finite set $x\subset T$ the tree $T\restriction
x$ is defined as the set of all elements of $T$ $\subset$-comparable with some
element of $x$
\item for a tree $S\in\nwd$ and a sequence $\eta\in S$ the tree $S(\eta)$ is 
the set $\{\tau:\eta^\smallfrown\tau\in S\}.$
\item $\um,$ the universal meager forcing \cite {S, Definition 4.2}, is the set
$\{\langle n,S\rangle:n\in\omega,S\in\nwd\}$ ordered by $\langle n,S\rangle
\leq \langle m,T\rangle$ if $n\geq m,T\subset S$ and $S\cap{^m2}=
T\cap {^m2}.$
\item for a tree $U\in\nwd$ we write $\um(U)=\{\langle n,S\rangle\in\um:
S\subset U\};$ this set has an order on it inherited from $\um$
\endroster
\enddefinition

Obviously, the colection $\{[T]:T\in\nwd\}$ is a base for the nowhere
dense ideal. $\um$ is a $\sigma$-centered Souslin forcing designed
so as to produce a very large nowhere dense tree: if $G\subset\um$
is a generic filter, then this tree is $U_G=\bigcup
\{T:\langle 0,T\rangle\in G\};$ it is nowhere dense and it codes
the generic filter. The following is the instrumental weakening of
genericity:

\definition {Definition 3.2}
Let $M$ be a transitive model of ZFC and $U\in\nwd.$ We say that
the tree $U$ $\um$-dominates the model $M$ if there is some
element $T\in M\cap\nwd$ included in $U$ and for every $T\in
M\cap\nwd$ there is an integer $n\in\omega$ such that setting
$x=T\cap U\cap {^n2},$ the inclusion $T\restriction x\subset
U\restriction x$ holds.
\enddefinition

This notion has certain obvious monotonicity properties.
Suppose $S\subset T$ and $U$ are perfect nowhere dense trees
and $n$ is an integer such that setting
$x=T\cap U\cap {^n2},$ $T\restriction x\subset
U\restriction x$ holds. Then with $y=S\cap U\cap{^n2}$ we have 
$S\restriction y\subset
U\restriction y$ and for any integer $m>n$ and
$z=T\cap U\cap {^m2}$ we have $T\restriction z\subset
U\restriction z.$

It is immediate that if $U\in\nwd$ is $\um$-generic then it $\um$-dominates
the ground model. On the other hand, any $\um$-dominating
tree covers an $\um$-generic tree:

\proclaim {Lemma 3.3}
Let $M$ be a transitive model of ZFC and let $U\in\nwd$ $\um$-dominate
the model $M.$ If $D\subset\um\cap M$ is a dense set in $M$ then 
$D\cap\um(U)$ is dense in $\um(U)\cap M.$
\endproclaim

\demo {Proof}
First, the set $\um(U)\cap M$ is nonempty. Now let $\langle n,S\rangle
\in\um(U)\cap M$ and let $D\in M$ be a dense subset of $\um\cap M$
which is an element of the model $M.$ We shall produce a condition
$p\in D\cap\um(U),p\leq \langle n,S\rangle,$ proving the lemma.

Work in $M.$ By induction on $i\in\omega,$ build conditions
$\langle n_i,T_i\rangle, p_{x,i}\in\um$ so that:

\roster
\item $\langle n_0,T_0\rangle=\langle n,S\rangle$ and 
$\langle n_{i+1},T_{i+1}\rangle\leq\langle n_i,T_i\rangle$
\item for every integer $i>0,$ for every sequence $\eta\in{^i2}$ there is
a sequence $\tau\in{^{n_i}2}$ with $\eta\subset\tau$ and $\tau\notin T_i$
\item to produce $\langle n_{i+1},T_{i+1}\rangle$ from
$\langle n_i,T_i\rangle,$ for every nonempty set $x\subset
{^{n_i}2}\cap T_i$ find a condition $p_{x,i}=\langle n_{x,i},S_{x,i}\rangle
\leq \langle n_i,T_i\restriction x_i\rangle$ in the dense set $D.$ Set $T_{i+1}=
\bigcup S_{x,i}$ and let $n_{i+1}$ be arbitrary so that (2) is satisfied.
\endroster

After this is done, let $T_\omega=\bigcup_i T_i.$ The induction hypothesis
(2) implies that $T_\omega\in\nwd\cap M$ and therefore, there is
an integer $i\in\omega$ such that setting $x={^{n_i}2}\cap T\cap U$
we get $T_\omega\restriction x\subset U\restriction x.$ Note that
the set $x$ is nonempty, because it includes $S\cap {^{n_i}2}.$
Now $p_{x,i}$ is the desired condition. \qed
\enddemo

\proclaim {Corollary 3.4}
(Subgenericity) Let $P$ be a forcing and $\dot S$ a $P$-name
such that

\roster
\item $P\Vdash$ the tree $\dot S\ \um$-dominates the ground model
\item for every $T\in\nwd$ the boolean value $||\check T\subset\dot S||_P$
is non-zero.
\endroster

Then there is a complete embedding $RO(\um)\lessdot P*(\um(\dot S)\cap$
the ground model$)=R$ such that $R\Vdash$``$\dot U\subset\dot
S$'', where $\dot U$ is the $\um$-generic tree.
\endproclaim

\proclaim {Corollary 3.5}
Let $M$ be a countable transitive model of ZFC such that $M\models$``$
P$ is a forcing adding a $\um$-dominating tree $\dot S$ as
in Corollary 3.4'', let $p\in P$ and let $T\in\nwd$ be any tree, 
not necessarily in the model
$M.$ Then there is an $M$-generic filter $G\subset P$ containing $p$ such that 
for some sequence $\eta\in\dot S/G$ we have $T\subset\dot S/G(\eta).$
\endproclaim

Set $\phi=$``cofinality of the meager ideal is $\aleph_1$''. The analysis
of the forcing $P_\phi$ is now completely parallel to the treatment
in Section 2. 

\definition {Definition 3.6}
A witness for $\phi$ is an $\omega_1$-sequence $s$ of perfect
nowhere dense trees such that

\roster
\item for every $\nwd$ tree $T$ the set $\{\alpha\in\omega_1:
T\subset s(\alpha)\}\subset\omega_1$ is stationary
\item for every $\nwd$ tree $T$ the set $C_T=\{\alpha\in\omega_1:$
there is $n\in\omega$ such that if $x={^n2}\cap T\cap s(\alpha)$
then $T\restriction x\subset s(\alpha)\restriction x\}\subset\omega_1$
contains a club
\item there is a $\nwd$ tree $T$ which is contained in all $s(\alpha):
\alpha\in\omega_1.$
\endroster
\enddefinition

Of course, it is important to verify that this notion deserves its name.

\proclaim {Lemma 3.7}
The following are equivalent:
\roster
\item the cofinality of the meager ideal is $\aleph_1$
\item the cofinality of the nowhere dense ideal is $\aleph_1$
\item there is a witness for $\phi.$
\endroster
\endproclaim

\demo {Proof}
(1)$\to$(2): let $\{Y_\alpha:\alpha\in\omega_1\}$
be a base for the meager ideal, $Y_\alpha\subset\bigcup\{[T_\alpha^i]:
i\in\omega\}$ for some sequence $T_\alpha^i:i\in\omega$ of $\nwd$
trees with $T_\alpha^i\subset T_\alpha^{i+1}.$ We shall show
that the collection $\{ [T_\alpha^i(\eta)]:\alpha\in\omega_1,
i\in\omega,\eta\in T_\alpha^i\}$ is a base for the nowhere dense ideal,
proving the lemma. Indeed, let $S$ be a nowhere dense tree on $\omega.$
We shall produce $\alpha,i,\eta$ so that $S\subset T_\alpha^i(\eta)$
and therefore $[S]\subset[T_\alpha^i(\eta)].$ 

It is a matter of an easy surgery on the tree $S$ to obtain a nowhere dense
tree $\bar S$ so that for every $\eta\in\bar S$ there is $\tau\in\bar S$
with $\eta\subset\tau$ and $\bar S(\tau)=S.$ Choose a countable ordinal
$\alpha$ so that $[\bar S]\subset Y_\alpha$ and attempt to build
 a descending sequence $\eta_i:i\in\omega$ of elements of $\bar S$
so that $\eta_i\notin T_\alpha^i.$ There must be an integer $i\in\omega$
such that the construction cannot proceed past $\eta_i$--otherwise
the branch $\bigcup_{i\in\omega}\eta_i\in[\bar S]$ would lie outside
of the set $Y_\alpha.$ But then, $\bar S(\eta_i)\subset T_\alpha^i(\eta_i)$
and if $\tau\in\bar S$ is such that $\eta_i\subset \tau$ and $\bar S(\tau)
=S$ then necessarily $S=\bar S(\tau)\subset T_\alpha^i(\tau).$ 

(2)$\to$(3): let $T_\xi:\xi\in\omega_1$  be a $\subset$-cofinal family
of $\nwd$ trees. Fix a partition $S_\xi:\xi\in\omega_1$ of $\omega_1$
into disjoint stationary sets and a $\nwd$ tree $T.$ For each
$\alpha\in\omega_1$ choose a filter $G_\alpha\subset\um$ such that

\roster
\item $\langle 0,T\cup T_\xi\rangle\in G_\alpha$ whenever $\alpha\in S_\xi$
\item $G_\alpha$ meets the dense sets $D_\beta=\{\langle n,U\rangle
\in\um:$ setting $x=U\cap {^n2}$ we have $T_\beta\restriction x\subset
U\restriction x\,$ for all $\beta\in\alpha$.
\endroster

Using the remarks after Definition 3.2 it is easy to see that the sequence
$s:\omega_1\to\nwd$ defined by $s(\alpha)=\bigcup\{ U:\langle 0,U\rangle\in
G_\alpha\}$ is the desired witness for $\phi.$

(3)$\to$(1) let $s:\omega_1\to\nwd$ be a witness for $\phi.$ Obviously,
the family $Y_\alpha=\bigcup_{\beta\in\alpha}[s(\beta)]:\alpha\in\omega_1$
is $\subset$-cofinal in the nowhere dense ideal. \qed 
\enddemo

Now if
$M$ is a transitive model of ZFC with $M\models$
``$s$ is a witness for $\phi$ and $\delta$ is a Woodin cardinal'' 
then in $M,$ $\Bbb Q_{<\delta}\Vdash$
``$js(\omega_1)$ is a $\nwd$ tree $\um$-covering the model $M$'',
where $j$ is the term for the generic nonstationary tower embedding;
also $M,\Bbb Q_{<\delta},js(\omega_1)$ 
satisfy the assumptions of Corollary 3.5. The proof of
$\Pi_2$-compactness of $\phi$ translates now literally from
the previous section. We prove the strategic iteration lemma from optimal
assumptions.

\proclaim {Strategic Iteration Lemma 3.8} Suppose that
the cofinality of the meager ideal is equal to $\aleph_1.$
The good player has a winning strategy in the game $\Cal G_\phi.$
\endproclaim

\demo {Proof}
Let $\vec N=\langle s,N_i,\delta_i:i\in\omega\rangle$ 
be a sequence of models with a witness for $\phi,$ let $y_0
\in\Bbb Q_{\vec N}$ and let $T$ be an arbitrary $\nwd$ tree.
We shall show that there is an $\vec N$-generic filter $G\subset\Bbb Q_{\vec N}$
such that letting $S=j_\Bbb Qs(\omega_1^{\vec N}),$ where $j_\Bbb Q$
is the $\Bbb Q_{<\delta_0}^{N_0}$-generic ultrapower embedding
using the filter $G\cap\Bbb Q_{<\delta_0}^{N_0},$ we have that for some
$\eta\in S,T\subset S(\eta).$ With this fact in hand, the winning strategy
for the good player consists just from an appropriate bookkeping:

Since cofinality of the meager ideal is $\aleph_1,$ it is possible to choose
a $\subset$-cofinal family $T_\alpha:\alpha\in\omega_1$ of $\nwd$
trees. So the good player can easily play the game so that with the
resulting embedding $j$ of the initial iterable model
$M,$ for every $\alpha\in\omega_1$ there is $\gamma
\in\omega_1$ and a sequence $\eta$ in the tree $
js(\gamma)$ such that $T_\alpha\subset
js(\gamma)(\eta).$ It is immediate that if this is the case and the iteration
$j$ is full, the sequence $j(s)$ is a witness for $\phi$ and the good
player won the run of the game. For let $S\in\nwd$ be an arbitrary tree.
Then there are $\alpha,\gamma$ and $\eta$ such that $S\subset T_\alpha
\subset js(\gamma)(\eta)$ and so

\roster
\item the set $\{\beta\in\omega_1:S\subset js(\beta)\}$ contains the set
$\{\beta\in\omega_1:js(\gamma)(\eta)\subset js(\beta)\}$ which is in the target
model of the iteration $j,$ is stationary there from Definition 3.6(1)
and so is stationary in $V$ by the fullness of the iteration
\item the set $\{\beta\in\omega_1:$ for some $n\in\omega,S\restriction x
\subset js(\beta)\restriction x$ holds with $x={^n2}\cap S\cap js(\beta)\}$
contains the set $\{\beta\in\omega_1:$ for some $n\in\omega,js(\gamma)(\eta)
\restriction x
\subset js(\beta)\restriction x$ holds with $x={^n2}\cap js(\gamma)(\eta)
\cap js(\beta)\},$ which is in the target model of the iteration
$j$ and contains a club by Definition 3.6(2).
\endroster

Therefore Definition 3.6(1,2) are verified for $j(s)$ and (3) of that
definition follows from elementarity of the embedding $j.$ Thus
$js$ is a witness for $\phi$ as desired.

The proof of the local fact about the sequence of models carries over from
Lemma 2.8 with the following changes:

\roster
\item $d$ is replaced with $s,$
$\Bbb D$ is replaced with $\um,$ the Hechler real $e$ is replaced
with a $\nwd$ tree $U$
\item the step (5) of that proof is replaced with: there is a sequence
$\eta\in U$ such that $T\subset U(\eta).$ Note that the set
$\{\langle n,S\rangle\in\um:\exists\eta\in S\ T\subset S(\eta)\}$
is dense in $\um.$
\item the ordering $R_i$ is defined as follows: $R_i=\{\langle
y,n,S\rangle:y\in\Bbb Q_i,\langle n,\ S\rangle\in\um,\ S=s(\alpha)\restriction
z$ for some $\alpha\in\omega_1^{\vec N}$ and finite set $z$ such that
$\forall x\in y\ S\subset s(x\cap\omega_1)\}.$ It is possible to restrict
ourselves to the trees $S$ of the above form, since the sequence $s(\alpha):
\alpha\in\omega_1$ is $\subset$-cofinal in $\nwd\cap N_i.$
\endroster
\qed
\enddemo

\proclaim {Conclusion 3.9}
The sentence $\phi=$ cofinality of the meager ideal$=\aleph_1$
is $\Pi_2$-compact. Theorem Scheme 0.2 holds even with an extra predicate for
$\omega_1$-sequences of meager sets cofinal in the ideal.
\endproclaim

\subhead
{3.1. Cofinality of the null ideal}
\endsubhead

In this subsection we shall show that ``cofinality of the null ideal
$=\aleph_1$'' is a $\Pi_2$-compact statement. The following textbook
equality will be used:

\proclaim {Lemma 3.10}
Cofinality of the null ideal is equal to the cofinality of the poset of
the open subsets of reals of finite measure ordered by inclusion.
\endproclaim

Therefore we will really care about large open sets of finite measure. 

\proclaim {Definition 3.11}
The amoeba forcing $\Bbb A$ is the set $\{\langle\Cal O,\epsilon\rangle:\Cal O$
is an open set of finite measure and $\epsilon$ is a positive rational
greater than $\mu(\Cal O)\}$ ordered by $\langle \Cal O,\epsilon\rangle
\leq\langle\Cal P,\delta\rangle$ if $\Cal P\subset\Cal O$ and $\epsilon
\leq\delta.$ The restricted poset $\Bbb A(\Cal O)$ for an open set
$O\subset\Bbb R$ is $\{\langle\Cal P,\epsilon\rangle\in A:\Cal P\subset
\Cal O\}$ with the inherited ordering.
\endproclaim

It is not a priori clear why the different versions of the amoeba forcing
should be isomorphic, see \cite {Tr}.
The amoeba poset is a $\sigma$-linked Souslin forcing notion designed
to add a ``large'' open set of finite measure. If $G\subset \Bbb A$ is generic
then the set $\Cal O_G=\bigcup\{\Cal P:\langle\Cal P,\epsilon\rangle
\in G$ for some $\epsilon\}$ is this open set and it determines the generic
filter. Again, there is a natural weakening of the notion of
genericity. Fix once and for all a sequence $f_i:i\in\omega$
of measure-preserving functions from $\Bbb R$ to $\Bbb R$ so that
the sets $f_i''\Bbb R$ are pairwise disjoint and the sequence is
arithmetical.

\proclaim {Definition 3.12}
Let $M$ be a transitive model of ZFC and $\Cal O$ be an open set
of reals. We say that $\Cal O$ $\Bbb A$-dominates the model $M$
if for every open set $\Cal P$ of finite measure in the model $M$
for all but finitely many integers $m\in\omega,$ $f_m^{-1}\Cal P\subset\Cal O.$
\endproclaim

Obviously, the amoeba generic open set does $\Bbb A$-dominate the ground
model. We aim for the subgenericity theorems.

\proclaim {Lemma 3.13}
Let $M$ be a transitive model of ZFC and let $\Cal O$ dominate $M.$
For every dense set $D\subset \Bbb A\cap M$ which is in the model $M$
the set $D\cap\Bbb  A(\Cal O)$ is dense in $M\cap\Bbb A(\Cal O).$
\endproclaim

\demo {Proof}
Let $M,\Cal O,D$ be as in the lemma and let $\langle\Cal P,\epsilon\rangle
\in M\cap\Bbb A(\Cal O).$ 
We shall produce a condition $p\in D\cap\Bbb A(\Cal O)$
below $\langle\Cal P,\epsilon\rangle,$ proving the lemma.

Work in the model $M.$ By induction on $i\in\omega$ build conditions
$\langle\Cal R_i,\delta_i\rangle\leq\langle \Cal P,\epsilon_i\rangle$
so that:

\roster
\item $\epsilon=\epsilon_0,\langle\Cal R_i,\delta_i\rangle\in D$
\item for every integer $i>0$ the inequality $\epsilon_i-\mu
(\Cal P)<2^{-i}$ holds.
\endroster

Let $\Cal S\in M$ be any open set of finite measure which covers
the set $\bigcup_{i\in\omega}f_i(\Cal R_i\setminus \Cal P).$ Since
$\Cal O$ $\Bbb A$-dominates the model $M,$ there must be an integer
$i\in\omega$ such that $f_i^{-1}\Cal S\subset\Cal O,$ and so
$\Cal R_i\subset\Cal O.$ Then
$\langle\Cal R_i,\delta_i\rangle\leq\langle \Cal P,\epsilon_i\rangle$
is the desired condition. \qed
\enddemo

\proclaim {Corollary 3.14}
(Subgenericity) Let $P$ be a forcing and $\dot\Cal O$ a
$P$-name such that

\roster
\item $P\Vdash$``the open set $\dot\Cal O\subset\Bbb R$ $\Bbb A$-dominates
the ground model''
\item for every open set $\Cal P$ of finite measure the boolean value
$||\check\Cal P\subset\dot\Cal O||_P$ is nonzero.
\endroster

Then there is a complete embedding $RO(\Bbb A)\lessdot P*(\Bbb A(\dot O)\cap$
the ground model$)=R$ such that $R\Vdash$``$\dot\Cal R\subset
\dot\Cal O$'', where $\dot\Cal R$ is the name for the $A$-generic
open set.
\endproclaim

A witness for $\phi=$``cofinality of the null ideal
$=\aleph_1$'' is an $\omega_1$-sequence $o$
of open sets of finite measure such that

\roster
\item for every open set $\Cal P$ of finite measure the set
$\{\alpha\in\omega_1:\Cal P\subset o(\alpha)\}\subset\omega_1$ is stationary
\item for every open set $\Cal P$ of finite measure the set
$\{\alpha\in\omega_1:$ for all but finitely many integers $m\in\omega$
$f_m^{-1}\Cal P\subset o(\alpha)\}$ contains a club.
\endroster

Again, it is very simple to prove using Lemma 3.10 that $\phi$ is equivalent
with the existence of a witness. The analysis of the forcing
$P_\phi$ almost literally translates from Section 2. We leave all
of this to the reader.

\proclaim {Conclusion 3.15}
The statement $\phi=$``cofinality of the null ideal is $\aleph_1$'' is
$\Pi_2$-compact. Theorem scheme 0.2 holds even with a predicate
for cofinal families of null sets added to the language of $\stru.$
\endproclaim

\head
{4. Souslin Trees}
\endhead

The assertion ``there is a Souslin tree'' does not seem to be $\Pi_2$-compact
as outlined in Subsection 1.3; however, some of its variations are.
A $P_{max}$-style model in which many Souslin trees exist was constructed
in \cite {W2} and in the course of the argument the following theorem,
which implies the strategic iteration lemmas for all
sentences considered in this section, was proved.

Let $\Cal G_S$ be a two-person game played along the lines of 
$\Cal G_\phi$--defined in 1.12--with the following modifications:

\roster
\item in the 0-th move the player Bad specifies a collection $\Cal S$
of Souslin trees in the model $M$ instead of just one witness
\item in the $\alpha$-th step the player Bad must choose a sequence
$\langle N_i:i\in\omega\rangle$ of models such that $j_{\gamma_\alpha}(\Cal S)$
consists of Souslin trees as seen from each $N_i:i\in\omega$
\item the player Good wins if $j_{\omega_1}(\Cal S)$ is a collection
of Souslin trees.
\endroster

\proclaim
{Strategic Iteration Lemma 4.1} \cite {W2}
Assume $\lozenge.$ Then the player Good
has a winning strategy in the game $\Cal G_S.$
\endproclaim

\demo {Proof}
Recall that $\omega_1$-trees are by our convention sets of functions
from countable ordinals to $\omega$ with some special properties.
Fix a diamond sequence $\langle A_\beta:\beta\in\omega_1\rangle$
guessing uncountable subsets of such trees.
The player Good wins the game as follows. Suppose we are at $\alpha$-th
stage of the play and let $\beta=\omega_1^{M_{\gamma_\alpha}}$ and
$\Cal S_\alpha=j_{0,\gamma_\alpha}(\Cal S).$ Suppose Bad played
a sequence $\vec N=\langle N_i,\delta_i:i\in\omega\rangle$ of models
according to the rules--so $N_0=M_{\gamma_\alpha}$--and some $p\in
\Bbb Q_{\vec N}.$ Let us call an $\vec N$-generic filter $G\subset
\Bbb Q_{\vec N}$ good if setting $j_{\Bbb Q}$ to be the ultrapower
embedding of $N_0$ 
derived from $G\cap\Bbb Q_0$ we have: for every tree $S\in\Cal S_\alpha,$
if $A_\beta\subset S$ is a maximal antichain then every node
at $\beta$-th level of $j_{\Bbb Q}(S)$ has an element of $A_\beta$
above it in the tree ordering.

If the player Good succeeds in playing a good filter containing $p$
at each stage $\alpha\in\omega_1$ of the game then he wins: every tree
in the collection $j_{0\omega_1}(\Cal S)$ can then be shown Souslin by
the usual diamond argument. Thus the following claim completes the proof.

\proclaim {Claim 4.2}
At stage $\alpha$ there is a good filter $G\subset\Bbb Q_{\vec N}$ containing
$p.$
\endproclaim

\demo {Proof}
Actually, any sufficiently generic filter is good. Note that every
$\Bbb Q_0$ name $\dot y\in N_0$ for a cofinal branch of any tree $S\in
\Cal S_\alpha$ is in fact a $\Bbb Q_{\vec N}$-name for a generic subset of 
$S$--this follows from the fact that $S$ is a Souslin tree in every model
$N_i:i\in\omega.$ Thus if a filter $G\subset\Bbb Q_{\vec N}$ meets
every dense set recursive in some fixed real coding $\vec N$ and $A_\beta,$
necessarily the branch $\dot y/G$ meets the set $A_\beta$ if $A_\beta\subset
S$ is a maximal antichain. Consequently, such a filter is good, since
every $\Bbb Q_0$ name $\dot y\in N_0$ for an element of $\beta$-th
level of $j_\Bbb Q(S)$ can be identified with a name for a cofinal
branch of the tree $S.$ \qed
\enddemo
\enddemo

\subhead
{4.0. Free Souslin trees}
\endsubhead

The first $\Pi_2$-compact sentence considered in this section is
$\phi=$``there is a free tree'' as clarified in the following definition:

\definition
{Definition 4.3} A Souslin tree $S$ is free if for every finite collection
$s_i:i\in I$ of distinct elements of the same level of $S$ the forcing
$\prod_{i\in I}S\restriction s_i$ is c.c.c.
\enddefinition

Thus every finitely many pairwise distinct cofinal branches of a free tree are
mutually generic.
It is not difficult to prove that both of the classical methods for forcing
a Souslin tree \cite {Te, J1} in fact provide free trees. It is an open problem
whether existence of Souslin trees implies existence of free trees.

The following observation, pointed out to us by W. Hugh Woodin, 
greatly simplifies
the proof of the optimal iteration lemma for $\phi:$ any sufficiently
rich (external) collection of cofinal branches of a free tree
determines a {\it symmetric extension} of the universe in the
appropriate sense.

\proclaim
{Lemma 4.4} Suppose that $M$ is a countable transitive model of a rich fragment
of ZFC, $M\models$``$S$ is a free Souslin tree'' and $B=\{b_i:i\in I\}$
is a countable collection of cofinal branches of $S$ such that $\bigcup
B=S.$ Then there is an enumeration $b_j:j\in\omega$ of $B$ such that
the equations $b_j=\dot c_j$ determine an $M$-generic filter on $P_S.$
\endproclaim

Here, $P_S\in M$ is the finite support product of countably many copies
of the tree $S,$ with $\dot c_j:j\in\omega$ being the canonical $P_S$-names
for the added $\omega$ branches of $S.$

\proclaim
{Corollary 4.5} Suppose $M,S$ and $B$ are as in the Lemma and suppose that
$M\models$``$P$ is a forcing, $p\in P$ and $P\Vdash\dot C$ is a collection
of cofinal branches of the tree $S$ such that $\bigcup\dot C=\check S$''.
Then there is an $M$-generic filter $G\subset P$ with $p\in G$ and
$\dot C/G=B.$
\endproclaim

\demo {Proof}
Work in $M.$ Without loss of generality 
we may assume that $p=1$ and that the forcing
$P$ collapses both $\kappa=(2^{\aleph_1})^+$ and $|\dot C|$ to $\aleph_0.$
(Otherwise switch to $P\times Coll(\omega,\lambda)$ for 
some large enough ordinal $\lambda$.)
There is a complete embedding of $RO(P_S)$ into $RO(P)$ such that
$P\Vdash$``$\dot C$ is the canonical set of branches of $S$ added by $P_S$
under this embedding''. This follows from Lemma 4.4 applied in $M^P$ to
$M\cap H_\kappa,S$ and $\dot C.$ Another application of the Lemma to
$M,S$ and $B$ gives an $M$-generic filter $H\subset P_S$ such that
$B$ is the canonical set of branches of the tree $S$ added by $H.$
Obviously, any $M$-generic filter $G\subset P$ extending $H$--via
the abovementioned embedding--is as desired. \qed
\enddemo

\demo
{Proof of Lemma} Say that the conditions in $P_S$ have the form of
functions from some $n\in\omega$ to $S$ with the natural ordering.
We shall show that for each injective $f:n\to B$ and every open dense
set $O\subset P_S$ in the model $M$ there is an injection $g:m\to B$
extending $f$ and a condition $p\in O$ with $\dom(p)=m$ and
$\bigwedge_{k\in m}p(k)\in g(k).$ Granted that, a construction
of the desired enumeration is straigthforward by the obvious bookkeeping
argument using the countability of both $M$ and $B.$

So fix $f$ and $O$ as above. There is an ordinal $\alpha\in\omega_1^M$
such that the branches $f(k):k\in n$ pick pairwise distinct elements
$s_k:k\in n$ from $\alpha$-th level of the tree $S.$ Let $D=
\{z\in\prod_{k\in n}S\restriction s_k:\exists p\in O\ p\restriction n=z\}
\in M.$ Since $O\subset P_S$ is dense below the condition $\langle s_k:
k\in n\rangle\in P_S,$ the set $D$ must be dense in 
$\prod_{k\in n}S\restriction s_k.$ Since this product is c.c.c. in
the model $M,$ the branches $f(k):k\in n$ determine an 
$M$-generic filter on it and there must be $z\in D$ such that
$\bigwedge_{k\in n}z(k)\in f(k).$ Choose a condition $p\in O$
with $\dom(p)=m$ and $p\restriction n=z.$ Since $\bigcup B=S,$ it is possible
to find branches $g(k):n\leq k<m$ in the set $B$ which are pairwise distinct
and do not occur on the list $f(k):k\in n$ such that $\bigwedge_{n\leq k<m}
p(k)\in g(k).$ The branches $f(k):k\in n$ and $g(k):n\leq k<m$ together
give the desired injection. \qed
\enddemo

\proclaim
{Optimal Iteration Lemma 4.6} Assume there is a free tree. Whenever
$M$ is a countable transitive model of ZFC iterable with
respect to its Woodin cardinal $\delta$ and $M\models$``$U$ is
a free tree'' there is a full iteration $j$ of $M$ such that
$j(U)$ is a free tree.
\endproclaim

\demo {Proof}
Let $T$ be a free Souslin tree and let $M,U,\delta$
be as above; so $M\models$``$U$ is a free Souslin tree."
We shall produce a full iteration $j$
of $M$ such that there is a club $C\subset\omega_1$ and an isomorphism
$\pi:T\restriction C\to j(U)\restriction C.$ Then, since the trees $T,j(U)$ are
isomorphic on a club, necessarily $j(U)$ is a free Souslin tree. This will
finish the proof of the lemma.

The iteration will be constructed by induction on $\alpha\in\omega_1$ and 
we will have $\theta_\alpha=\omega_1^{M_\alpha}$ and $C=\{\theta_\alpha:
\alpha\in\omega_1\}.$  Also, we shall write $U_\alpha$ for
the image of the tree $U$ under the embedding $j_{0,\alpha}.$
This is not to be confused with the $\alpha$-th level of the
tree $U.$ In this proof, levels of trees are never indexed
by the letter $\alpha.$ 

First, fix a
partition $\{S_\xi:\xi\in\omega_1\}$ of the set of countable limit ordinals
into disjoint stationary sets. By induction on $\alpha\in\omega_1,$
build the models together with the elementary embeddings,
plus an isomorphism $\pi:T\restriction C\to jU\restriction C,$ plus
an enumeration $\{\langle x_\xi,\beta_\xi\rangle:\xi\in\omega_1\}$ 
of all pairs $\langle x,\beta\rangle$ with $x\in\Bbb Q_\beta.$
The induction hypotheses at $\alpha\in\omega_1$ are:

\roster
\item the function $\pi\restriction T\restriction \{\theta_\gamma:
\gamma\in\alpha\}$ has been defined and it is an isomorphism of
$T\restriction \{\theta_\gamma:\gamma\in\alpha\}$ and $U_\alpha
\restriction \{\theta_\gamma:\gamma\in\alpha\}$ 
\item the initial segment $\{\langle x_\xi,\beta_\xi\rangle:
\xi\in\theta_\alpha\}$ has been constructed and every pair
$\langle x,\beta\rangle$ with $\beta\in\alpha$ and $x\in\Bbb Q_\beta$
appears on it
\item if $\gamma\in\alpha$ belongs to some--unique--set $S_\xi$
then $j_{\beta_\xi,\gamma}(x_\xi)\in G_\gamma.$
\endroster

At limit steps, we just take direct limits and unions. At successor steps,
given $M_\alpha,U_\alpha,\delta_\alpha$ and $\pi\restriction\{
\theta_\gamma:\gamma\in\alpha\},$ we must provide an $M_\alpha$-generic
filter $G_\alpha\subset\Bbb Q_\alpha$ such that
setting $U_{\alpha+1}=j_{\Bbb Q}U_\alpha,$ 
where $j_{\Bbb Q}$ is the generic ultrapower
of $M_\alpha$ by $G_\alpha,$ it is possible to extend the isomorphism
$\pi$ to $\theta_\alpha$-th levels of $T$ and $U_{\alpha+1}.$

 First suppose $\alpha$ is a successor ordinal, $\alpha=\beta+1.$
Let $G_\alpha$ be an arbitrary $M_\alpha$-generic filter; we claim that
$G_\alpha$ works. Simply let  for every $t\in T_{\theta_\beta} $
$\pi\restriction(T\restriction t)_{\theta_\alpha}$ to be a bijection
of $(T\restriction t)_{\theta_\alpha}$ and $(U_{\alpha+1}\restriction
\pi(t))_{\theta_\alpha}.$ This is clearly possible since
both of these sets are infinite and countable. Induction hypothesis
(1) continues to hold, induction hypothesis (2) is easily arranged by
extending the enumeration properly and (3) does not say anything
about successor ordinals.

Finally, suppose $\alpha$ is a limit ordinal. In this case, the $\theta_
\alpha$-th level of the tree $U_{\alpha+1}$ is determined by
$\pi\restriction T\restriction\{\theta_\gamma:\gamma\in\alpha\}$
and the necessity of extending the isomorphism $\pi$ to the $\theta_
\alpha$-th level of the tree $T.$ Namely we must have
$\theta_\alpha$-th level of $U_{\alpha+1}$ equal to the set
$D=\{d_t:t\in T_{\theta_\alpha}\}$ where $d_t=\bigcup\{
\pi(r):r\in T\restriction \{\theta_\gamma:\gamma\in\alpha\},
t\leq_Tr\}.$ Corollary 4.5 applied to $M_\alpha,S_\alpha,\Bbb Q_\alpha,
(S_{\alpha+1})_{\theta_\alpha}$ and $D$
shows that an appropriate generic filter on $\Bbb Q_\alpha$ can be found
containing the condition 
$j_{\beta_\xi,\alpha}(x_\xi)$ if $\alpha\in S_\xi.$ The isomorphism
$\pi$ then extends in the obvious unique fashion mapping $t$ to $d_t.$ \qed
\enddemo 

\proclaim
{Conclusion 4.7}
 The sentence $\phi=$``there is a free tree'' is $\Pi_2$-compact.
\endproclaim

Another corollary to the proof of Lemma 4.6 is the fact that $\Sigma_1^1$
theory of free trees is {\it complete} and {\it minimal} in the
following sense. Suppose $\psi$ is a $\Sigma_1^1$ property
of trees $T$ which depends only on the Boolean algebra $RO(T),$
that is, ZFC$\vdash$``$RO(S)=RO(T)$ implies $T\models\psi$ iff $S\models\psi$
''. Then, granted large cardinals, the sentence $\psi$ is either
true on all free trees in all set-generic extensions of the universe
or it fails on all such trees. Moreover, if $\psi$ fails on
{\it any} $\omega_1$-tree in any set-generic extension then it fails
on all free trees. Here, by $\Sigma_1^1$ property we mean
a formula of the form $\exists A\subset T\chi,$ where all quantifiers
of $\chi$ range over the elements of $T$ only.

It should be noted that it is impossible to add a predicate $\frak S$
for free trees to the language of $\stru$ and preserve the compactness
result. For consider the following two $\Pi_2$ sentences for
$\langle H_{\aleph_2},\in,\frak I,\frak S\rangle:$

$$\eqalign{\psi_0&=\forall S,T\in\frak S\ \exists s\in S,t\in T\ S\restriction s\times
T\restriction t\text{ is c.c.c.}\cr
\psi_1&=\text{for every Aronszajn tree }T\text{ there is a tree }S\in\frak S
\text{ such that }T\Vdash S\text{ is special.}\cr}$$

It is immediate that $\psi_0$ and $\psi_1$ together imply that $\frak S$ is empty, i.e. $\lnot\phi$.
Meanwhile,
$\psi_0\land\phi$ was found consistent in \cite {AS}--and in fact holds in
our model--and $\psi_1\land\phi$ holds after adding $\aleph_2$ Cohen
reals to any model of GCH, owing to the following lemma:
 
\proclaim {Lemma 4.8}
For every aronszajn tree $T,$ $\Bbb C_{\aleph_1}\Vdash$``there is a free tree which is
specialized after forcing with $T$''.
\endproclaim

Note that Cohen algebras preserve Souslin trees.

\demo {Proof}
Let $T$ be an Aronszajn tree. Define a forcing $P$ as a set of pairs
$p=\langle s_p,f_p\rangle$ where

\roster
\item $s_p$ is a finite tree on $\omega_1\times\omega$ such that $\langle
\alpha,n\rangle<_{s_p}\langle\beta,m\rangle$ implies $\beta<\alpha$ ...
this is a finite piece of the tree $S$ under construction
\item $f_p$ is a finite function with domain contained in $T$ and each
$f_p(t)$ a function from $\dom(s)\cap(\{\alpha\}\times\omega)$ to $\omega$
where $\alpha=lev(t)$ ...this is a finite piece of the $S$-specializing
$T$-name
\item for every $i<_{s_p}j$ and $t<_Tu$ the inequality $f_p(t)(i)\neq
f_p(u)(j)$ holds, if the relevant terms are defined ...this is the
specializing condition.
\endroster

The ordering is defined by $q\leq p$ if $\dom(s_p)\subset\dom(s_q)$ and
$s_q\cap\dom(s_p)\times\dom(s_p)=s_p$ and $f_p(t)\subset f_q(t)$
whenever $t\in\dom(f_p).$

Let $G\subset P$ be a generic filter and in $V[G]$ define a tree $S$
on $\omega_1\times\omega$ as the unique tree extending all
$s_p:p\in G,$ and a function $\tau$ on the tree $T$ to be $\tau(t)
=\bigcup_{p\in G}f_p(t).$ Obviously, $\tau$ represents a $T$-name
for a specializing function on $S:$ if $b\subset T$ is a cofinal branch
then the function $g:S\to\omega,$ $g=\bigcup _{t\in b}\tau(t)$
specializes the tree $S$ due to the condition (3) in the definition of
the forcing $P.$ To complete the proof, we have to verify that $RO(P)=
\Bbb C_{\aleph_1}$ and that $P\Vdash S$ is a free Souslin tree. This is done
in the following two claims.

\proclaim {Claim 4.9}
$RO(P)$ is isomorphic to $\Bbb C_{\aleph_1}.$
\endproclaim

\demo {Proof}
Obviously, $P$ has uniform density $\aleph_1,$ therefore it is enough to
prove that $P$ has a closed unbounded collection of regular subposets
\cite {K}. Let $\alpha\in\omega_1$ be a limit ordinal and let
$P_\alpha=\{ p\in P:dom(s_p)\subset\alpha\times\omega\}.$ It is easy to
verify that all $P_\alpha$'s are regular subposets of $P$ and that
they constitute an increasing continuous chain exhausting all of $P,$
proving the lemma. \qed
\enddemo

\proclaim {Claim 4.10}
$P\Vdash S$ is a free Souslin tree.
\endproclaim

\demo {Proof}
Assume that $p\Vdash$``$\dot A=\{a_\alpha:\alpha\in\omega_1\}$ is a family
of pairwise distinct elements of $\dot S\restriction i_0\times
\dot S\restriction i_1\times\dots\times\dot S\restriction i_n$'',
for some integer $n$ and pairwise $s_p$-incompatible elements $i_0\dots i_n$
of $\dom(s_p).$ To prove the lemma, it is enough to produce a condition
$q\leq p$ and ordinals $\alpha<\beta$ such that $q\Vdash \dot a_\alpha$ and
$\dot a_\beta$ are compatible.

Pick $p_\alpha,a_\alpha:\alpha\in\omega_1$ such that each $p_\alpha$
is a condition stronger than $p$ and it decides the value of the
name $\dot a_\alpha$ to be $a_\alpha,$ regarded as an $n+1$-element
subset of $\dom(s_{p_\alpha})\restriction i_0\cup\dom(s_{p_\alpha})
\restriction i_1\cup\dots\cup
   \dom(s_{p_\alpha})\restriction i_n.$

By a repeated use of counting arguments and a $\Delta$-system lemma,
thinning out the collection of $p_\alpha,a_\alpha$'s we may assume that

\roster
\item $\dom(s_{p_\alpha})$ form a $\Delta$-system with root $r$ and
$s_{p_\alpha}\restriction r\times r$ is the same for all $\alpha$
\item even the sets $lev(s_{p_\alpha})=
\{\beta\in\omega_1:\dom(s_{p_\alpha})\cap\{\beta\}\times\omega\neq 0\}$
form a $\Delta$ system with root
$lev(r)=\{\beta\in\omega_1:r\cap\{\beta\}\times\omega\neq 0\}.$
\item $f_{p_\alpha}\restriction T_\beta$ are the same
for all $\alpha,$ this for all $\beta\in lev(r)$
\item $a_\alpha$ form a $\Delta$-system with root $b\subset r.$
\endroster

Now let $x_\alpha=
\dom(f_{p_\alpha})\setminus \bigcup_{\beta\in lev(r)} T_\beta.$
Thus $x_\alpha$ are pairwise disjoint finite subsets of the Aronszajn
tree $T,$ and it is possible to find  countable ordinals $\alpha<\beta$
such that every $t\in x_\alpha$ is $T$-incompatible with every $u\in x_\beta.$
It follows that any tree $s_q$ with $\dom(s_q)=\dom(s_{p_\alpha})\cup
\dom(s_{p_\beta}),$ $s_q\restriction\dom(s_{p_\alpha})\times\dom(s_{p_\alpha})=
s_{p_\alpha}$ and $s_q\restriction\dom(s_{p_\beta})\times\dom(s_{p_\beta})=
s_{p_\beta},$ together with the function $f_q=f_{p_\alpha}\cup f_{p_\beta}$
give a condition $q=\langle s_q,f_q\rangle$ in the forcing $P$ which is
stronger than both $p_\alpha$ and $p_\beta.$ It is a matter of an easy surgery
on $s_{p_\alpha}$ and $s_{p_\beta}$ to provide such a condition
$q$ so that $a_\alpha,a_\beta$ are compatible in $s_q\restriction
i_0\times s_q\restriction i_1\times\dots\times s_q\restriction i_n.$
Then $p\geq q\Vdash$``$\dot a_\alpha$ and $\dot a_\beta$ are compatible
elements of $\dot A$'' as desired. \qed
\enddemo 

\enddemo
\subhead
{4.1. Strongly homogeneous Souslin trees}
\endsubhead

In this subsection it is proved that the assertion $\phi=$``there is a 
strongly homogeneous Souslin tree'' is $\Pi_2$-compact, where

\definition {Definition 4.11}
Let $T$ be an $\omega_1$-tree. A family $\{ h(s_0,s_1):s_0,s_1\in T$
are elements of the same level of $T\}$ is called
{\it coherent} if:
\roster
\item $h(s_0,s_1)$ is a level-preserving isomorphism of $T\restriction s_0$
and $T\restriction s_1;$ $h(s,s)=id$
\item (commutativity) let $s_0,s_1,s_2$ be elements of the same level
of $T$ and $t_0\leq s_0.$ Then $h(s_1,s_2)h(s_0,s_1)(t_0)=
h(s_0,s_2)(t_0)$
\item (coherence) let $s_0,s_1$ be elements of the same level of $T$ and
$t_0\leq s_0.$ Let $t_1=h(s_0,s_1)(t_0)\leq s_1.$ Then $h(t_0,t_1)=
h(s_0,s_1)\restriction T\restriction t_0$
\item (transitivity) if $\alpha$ is a limit ordinal and $t_0,t_1$
are two different elements at $\alpha$-th level of $T$ then there are 
$s_0,s_1\in T_{<\alpha}$ such that $h(s_0,s_1)(t_0)=t_1.$
\endroster
A tree is called {\it strongly homogeneous} if it has a coherent family
of isomorphisms.
\enddefinition

The existence of strongly homogeneous Souslin trees can be proved from
$\lozenge$ by a standard argument. Also, Todorcevic's term for
a Souslin tree in one Cohen real extension provides in fact for
a strongly homogeneous tree:

\proclaim {Theorem 4.12}
$\Bbb C\Vdash$ there is a strongly homogeneous Souslin tree.
\endproclaim

\demo {Proof}
An elaboration on Todorcevic's argument \cite {T1}. Let $T$ be a family
of functions such that

\roster
\item every $f\in T$ is of the form $f:\alpha\to\omega,$ finite-to-one
for some countable ordinal $\alpha$
\item  for each $\alpha\in\omega_1$ there is $f\in T$ with $\alpha=
\dom(f)$
\item every two functions $f,g\in T$ are modulo finite equal on the
intersection of their domains
\item $T$ is closed under finite changes of its elements.
\endroster

Such a family is built as in \cite {T1} and it can be understood as a
special Aronszajn tree under the reverse inclusion order. If $c\in\rr$
is a Cohen real then \cite {T1} the tree $T_c=\{c\circ f:f\in T\}$
ordered by reverse inclusion is a Souslin tree in $V[c].$
To conclude the proof of the Theorem, we shall find a coherent family 
of isomorphisms of the tree $T$ which is easily seen to lift to the tree $T_c.$
Namely, let $f,g\in T,\dom(f)=\dom(g).$ Define $h(f,g)(e)=g\cup (e\setminus f)$
for $e\in T$ with $f\subset e.$ By (3) and (4) above this is a well-defined
function and an isomorphism of the trees $T\restriction f$ and $T\restriction
g.$ The easy proof that these isomorphisms form a coherent family on a tree 
$T$ which lifts to the tree $T_c$ is left to the reader. \qed
\enddemo

Paul Larson proved that every
strongly homogeneous Souslin tree contains a regularly embedded free tree. In fact, every strongly homogeneous
Souslin tree can be written as a product of two free trees. 

\proclaim 
{Optimal Iteration Lemma 4.13} Assume there is a strongly homogeneous
Souslin tree. Whenever $M$ is a countable transitive model of ZFC
iterable with respect to its Woodin cardinal $\delta$
with $M\models$``$U$ is a strongly homogeneous Souslin tree''
there is a full iteration $j$ of $M$ such that $j(U)$ is a
strongly homogeneous Souslin tree.
\endproclaim

\demo {Proof}
Let $T$ be a strongly homogeneous Souslin tree with a coherent family
$\{ g(t_0)(t_1):t_0,t_1\in T_\alpha$ for some $\alpha\in\omega_1\}$
of isomorphisms and let $M,U,\delta$
be as above and $M\models$``$U$ is a strongly homogeneous Souslin tree
as witnessed by a family $h=\{h(s_0)(s_1):s_0,s_1\in U_\xi$ for
some $\xi\in\omega_1^M\}".$ We shall produce a full iteration $j$
of $M$ such that there is a club $C\subset\omega_1$ and an isomorphism
$\pi:T\restriction C\to j(U)\restriction C$ which commutes with the internal
isomorphisms of the trees: $\pi g(t_0,t_1)(u)=jh(\pi t_0,\pi t_1)(\pi u)$ 
whenever
the relevant terms are defined. Then, since the trees $T,j(U)$ are
isomorphic on a club, necessarily $j(U)$ is a Souslin tree, and it is
strongly homogeneous as witnessed by the family $j(h).$ This will
finish the proof of the lemma. Again, below $U_\alpha$ denotes the tree $j_{0\alpha}U$ and not the $\alpha$-th level of $U$.
Levels of trees are never indexed by $\alpha.$

The iteration will be constructed by induction on $\alpha\in\omega_1$ and 
we will have $\theta_\alpha=\omega_1^{M_\alpha}$ and $C=\{\theta_\alpha:
\alpha\in\omega_1\}.$ First, fix a
partition $\{S_\xi:\xi\in\omega_1\}$ of the set of countable limit ordinals
into disjoint stationary sets. By induction on $\alpha\in\omega_1,$
we shall build the models together with the elementary embeddings,
plus an isomorphism $\pi:T\restriction C\to jU\restriction C,$ plus
an enumeration $\{\langle x_\xi,\beta_\xi\rangle:\xi\in\omega_1\}$ 
of all pairs $\langle x,\beta\rangle$ with $x\in\Bbb Q_\beta.$
The induction hypotheses at $\alpha\in\omega_1$ are:

\roster
\item the function $i\restriction T\restriction \{\theta_\gamma:
\gamma\in\alpha\}$ has been defined, it is an isomorphism of
$T\restriction \{\theta_\gamma:\gamma\in\alpha\}$ and $U_\alpha
\restriction \{\theta_\gamma:\gamma\in\alpha\}$ and it commutes
with the internal isomorphisms of the trees, i.e. $\pi g(t_0,t_1)(u)=
j_{0,\alpha}h(\pi t_0,\pi t_1)(\pi u)$ whenever the relevant terms are defined
\item the initial segment $\{\langle x_\xi,\beta_\xi\rangle:
\xi\in\theta_\alpha\}$ has been constructed and every pair
$\langle x,\beta\rangle$ with $\beta\in\alpha$ and 
$x\in\Bbb Q_\beta$
\item if $\gamma\in\alpha$ belongs to some--unique--set $S_\xi$
then $j_{\beta_\xi,\gamma}(x_\xi)\in G_\gamma.$
\endroster

As before, (1) is the crucial condition ensuring that the tree
$T$ is copied to $j(U)$ properly. (2,3) are just bookkeping requirements
for making the resulting iteration full.

At limit steps, we just take direct limits and unions of the isomorphisms and
enumerations constructed so far. At successor steps, given $M_\alpha,$
we must produce a $M_\alpha$-generic filter 
$G_\alpha\subset\Bbb Q_\alpha$
such that setting $\langle M_{\alpha+1},U_{\alpha+1},\delta_{\alpha+1}\rangle$
to be the generic ultrapower of $\langle M_\alpha,U_\alpha,\delta_\alpha
\rangle$ by $G_\alpha,$ the isomorphism $\pi$ can be extended to 
$\theta_\alpha$-th level of the trees $T$ and $U_{\alpha+1}$
preserving the induction hypothesis (1).

$\bullet${\it Case 1.} $\alpha$ is a successor ordinal, $\alpha=\beta+1.$ Choose
an arbitrary $M_\alpha$-generic filter 
$G_\alpha\subset\Bbb Q_\alpha.$ We shall show how the
isomorphism $\pi$ can be extended to the
$\theta_\alpha$-th level of the trees $T$ and $U_{\alpha+1}$
preserving the induction hypothesis (1).

Let $t\in T_{\theta_\alpha}$
be arbitrary. The $\theta_\beta$-{\it orbit} of $t$ is the set
$\{u\in T_{\theta_\alpha}:\exists t_0,t_1\in T_{\theta_\beta}\ 
u=g(t_0,t_1)(t)\}.$ By the commutativity property of the
isomorphisms $g,$ the $\theta_\alpha$-th level of the tree $T$
partitions into countably many disjoint $\gamma$-orbits $O_k:
k\in\omega.$ Also, for every $u\in T_{\theta_\beta}$ and integer
$k\in\omega$ there is a unique $t\in O_k$ with $t\leq_Tu.$
The same analysis applies to the tree $U_{\alpha+1}$ and isomorphisms $h.$
The $\theta_\alpha$-th level of the tree $U_{\alpha+1}$
partitions into countably many disjoint $\gamma$-orbits $N_k:
k\in\omega.$

Now  it is easy to see that there is a unique way 
to extend the function $\pi$ to
$T_{\theta_\alpha}$ so that $\pi''O_k=N_k$ and $\pi$ is
order-preseving. Such an extended function will satisfy the induction
hypothesis (1). The induction hypothesis (2) is easily
managed and the induction hypothesis (3) does not say anything about
successor ordinals $\alpha.$

$\bullet${\it Case 2.} 
$\alpha$ is a limit ordinal. In this case, $\theta_\alpha$-th
level of the tree $U_{\alpha+1}$ is already pre-determined by
$\pi\restriction T\restriction\{\theta_\gamma:\gamma\in\alpha\}$
and the necessity of extending $\pi.$ Namely, we must have
$(U_{\alpha+1})_{\theta_\alpha}=\{ u:$ there is $t\in T_{\theta_\alpha}$
such that $u=\bigcup\{\pi r:t\leq_T r\}\}.$ The challenge
is to find an $M_\alpha$-generic filter $G_\alpha\subset\Bbb Q_\alpha$
such that $(\dot U_{\alpha+1})_{\theta_\alpha}/G_\alpha$
is of the abovedescribed form. Then necessarily the only possible
orderpreserving extension of $\pi$ to $T_{\theta_\alpha}$ will
satisfy the induction hypothesis (1). We shall use the fact that
it is enough to know one element of $(U_{\alpha+1})_{\theta_\alpha}$
in order to determine the whole level--by transitivity, Definition 4.2(4).

Work in $M_\alpha.$ Fix $\dot u,$ an arbitrary 
$\Bbb Q_\alpha$-name for an element of 
$(U_{\alpha+1})_{\theta_\alpha},$ which will be identified with the
cofinal--and therefore $M_\alpha$-generic--branch of the tree $U_\alpha$
it determines. Let $b_0\in\Bbb Q_\alpha$ be defined as
$j_{\beta_\xi,\alpha}(x_\xi)$ if $\alpha$ belongs to some--unique--
set $S_\xi$ with $\xi\in\theta_\alpha,$ otherwise let $b_0=1$
in $\Bbb Q_\alpha$ Let $B$ be the complete
subalgebra of $RO(\Bbb Q_\alpha)$ generated by the name
$\dot u.$ By some Boolean algebra theory, there must be $b_1\leq b_0$
and $s\in U_\alpha$ so that $pr_{B}b_1=[[\check s\in\dot u]]_
{B}=b_2$ and $B\restriction b_2$ is isomorphic
to $RO(U_\alpha\restriction s)$ by an isomorphism generated
by the name $\dot u.$ Without loss of generality $lev(s)=\theta
_\gamma$ for some $\gamma\in\alpha,$ since the set $\{\theta_\gamma:
\gamma\in\alpha\}$ is cofinal in $\theta_\alpha.$

Now pick an element $t\in T_{\theta_\alpha}$ such that $t\leq_T\pi^{-1}s.$
Then $c=\bigcup\{\pi(r):t\leq_Tr\}$ is a cofinal $M_\alpha$-generic
branch through $U_\alpha$ containing $s.$ Let $H\subset B$ be the
$M_\alpha$-generic filter determined by the equation $c=\dot u$ and let
$G_\alpha\subset RO(\Bbb Q_\alpha)$ be any $M_\alpha$
generic filter with $H\subset G_\alpha,b_1\in G_\alpha.$ We claim that
$G_\alpha$ works.

Let $\langle M_{\alpha+1},U_{\alpha+1},\delta_{\alpha+1}\rangle$ be the
generic ultrapower of $\langle M_\alpha,U_\alpha,\delta_\alpha\rangle$
using the filter $G_\alpha.$ Define $\pi\restriction T_{\theta_\alpha}$
by $\pi g(r_0,r_1)(t)=j_{0,\alpha+1}h(s_0,s_1)(c)$ where
$r_0,r_1\in T_{\theta_\gamma}$ for some $\gamma\in\alpha$ and $t\leq_T r_0,$
$\pi(r_0)=s_0,\pi(r_1)=s_1$ and $t$ is the element of 
$T_{\theta_\alpha}$ used to generate
$c$ in the previous paragraph. By the induction hypothesis (1)
and coherence--Definition 4.2(4)--$\pi$ is well-defined, and by transitivity
applied to both $T$ and $U$ side $\pi\restriction T_{\theta_\alpha}:
T_{\theta_\alpha}\to(U_{\alpha+1})_{\theta_\alpha}$ is a bijection.
It is now readily checked that $\pi$ commutes with the internal isomorphisms
$g,h$ and the induction hypothesis continues to hold at $\alpha+1.$
The hypothesis (2) is easily managed by suitably prolonging the enumeration,
and the induction hypothesis (3) is maintained by the choice of $b_0\in
G_\alpha.$ \qed
\enddemo
  
\proclaim {Conclusion 4.14}
The sentence ``there is a strongly homogeneous Souslin tree'' is
$\Pi_2$-compact.
\endproclaim
Again, the proof of Lemma 4.13 shows that the $\Sigma_1^1$ theory of
strongly homogeneous Souslin trees is complete in the same sense as 
explained in the previous Subsection. The first order theory of the
model obtained in this Subsection has been independently studied by Paul Larson.

\subhead
{4.2. Other types of Souslin trees}
\endsubhead

One can think of a great number of $\Sigma_1$ constraints on Souslin trees.
With each of them, the first two iteration lemmas can be proved for
$\phi=$``there is a Souslin tree with the given constraint'' owing to
Lemma 4.1. However, the absoluteness properties of the resulting models
as well as the status of $\Pi_2$-compactness of such sentences $\phi$
are unknown. Example:

\definition
{Definition 4.15}
A Souslin tree $T$ is self-specializing if $T\Vdash$``$\check T\setminus \dot b$
is special, where $\dot b$ is the generic branch''.
\enddefinition

A selfspecializing tree can be found under $\lozenge$ or after adding $\aleph_1$
Cohen reals. Such a tree is obviously neither free nor strongly
homogeneous and no such a tree exists in the models from the previous
two subsections.

\head
{5. The bounding number}
\endhead

In this section it will be proved that the sentence $\bb=\aleph_1$
is $\Pi_2$-compact even in the stronger sense with a predicate for
unbounded sequences added to the language of $\stru.$
Everywhere below, by an {\it unbounded sequence} we mean a modulo
finite increasing $\omega_1$-sequence of increasing functions
in $\rr$ without an upper bound in the eventual domination
ordering of $\rr.$

\subhead
{5.0. Combinatorics of $\bb$}
\endsubhead

A subgenericity theorem similar to the one obtained in the dominating number
section can be proved here too, in this case essentially saying that adding a
Cohen real is an optimal way of adding an unbounded real. It is just
a restatement of a familiar fact from recursion theory 
and is of limited use in what follows.

We abuse the notation a little writing $\Bbb C={^{<\omega}\omega}$ ordered by
reverse extension and for a function $f\in\rr,$
$\Bbb C(f)=\{\eta\in\Bbb C:\eta$ is on its domain
pointwise $\leq f\}$ ordered by reverse extension as well.

\proclaim {Lemma 5.1}
Let $M$ be a transitive model of ZFC and $f\in\rr$ be an increasing function
which is not bounded by any function in $M.$ Whenever $D\in M,$
$D\subset\Bbb C$ is a dense set, then $D\cap\Bbb C(f)\subset\Bbb C(f)$
is dense.
\endproclaim

\demo {Proof}
Fix a dense subset $D\in M$ of $\Bbb C$ and a condition $p\in \Bbb C(f).$
We shall produce $q\leq p,q\in D\cap\Bbb C(f),$ proving the lemma.

Work in the model $M.$ By induction on $n\in\omega$ build a sequence
$p_0,p_1\dots,p_n,\dots$ of conditions in $\Bbb C$ so that

\roster
\item $p_0=p,\dom(p_{n+1})>\dom(p_n)$
\item $p_{n+1}$ is any element of the set $D$ below the condition $r=
p^\smallfrown\langle 0,0\dots 0\rangle$ with as many zeros as necessary to get
$\dom(r)=m_n.$
\endroster

Let $X\subset\omega$ be the set $\{\dom(p_n):n\in\omega\}$ and let $g:
X\to\omega$ be given by $g(m_n)=\max(\rng(p_{n+1})).$ Then $X,g\in M$
and since the function $f$ is increasing and unbounded over the model $M,$
there must be an integer $n$ such that $f(\dom(p_n))>g(\dom(p_n)).$
Then obviously $q=p_{n+1}\leq p$ is the desired condition in
$D\cap\Bbb C(f).$ \qed
\enddemo

\proclaim {Corollary 5.2}
(Subgenericity) Let $P$ be a forcing and $\dot f$ a $P$-name such that

\roster
\item $P\Vdash$``$\dot f\in\rr$ is an increasing unbounded function''
\item for every finite sequence $\eta$ of integers the boolean value
$||\check \eta$ is bounded on its domain by $\dot f||_P$ is nonzero.
\endroster

Then there is a $P$-name $\dot Q$ and a complete embedding $\Bbb C\lessdot
RO(P*\dot Q)$ such that $P*\dot Q\Vdash$``$\dot f$ pointwise dominates
$\dot c,$ the $\Bbb C$-generic function''.
\endproclaim

\demo {Proof}
Set $\dot Q=\Bbb C(\dot f)$ and use the previous Lemma. \qed
\enddemo

It follows that a collection $A\subset\rr$ of increasing functions is unbounded
just in case the set $X=\{f\in\rr:$ some $g\in A$ eventually dominates
the function $f\}$ is nonmeager. For if $A$ is bounded by some $h\in\rr$
then $X\subset\{f\in\rr:h$ eventually dominates $f\}$ and the latter set
is meager; on the other hand, if $A$ is unbounded then Lemma 5.1 provides
sufficiently strong Cohen reals in the set $X$ to prove its nonmeagerness.
A posteriori, a forcing preserving nonmeager sets preserves unbounded
sequences as well.

A more important feature of the bounding number is that every two
unbounded sequences can be made in some sense isomorphic. Recall the
quasiordering $\leq_\bb$ defined in Subsection 1.2.

\definition
{Definition 5.3} For $b,c\in H_{\aleph_2}$ set $b\leq_\bb c$ if in every
forcing extension of the universe $b$ is an unbounded sequence implies
$c$ is an unbounded sequence.
\enddefinition

While under suitable assumptions (for example the Continuum Hypothesis)
the behavior of this quasiorder is very complicated, in the model
for $\bb=\aleph_1$ we will eventually build there will be exactly
two classes of $\leq_\bb$-equivalence. The key point is
the introduction of the following $\Sigma_1\stru$ concept
to ensure $\leq_\bb$-equivalence of two unbounded sequences.

\definition
{Definition 5.4} Unbounded sequences $b,c$ are locked if there is an infinite
set $x\subset\omega$ such that for every $\alpha\in\omega_1$ there is
$\beta\in\omega_1$ with $b(\beta)\restriction x$ eventually dominating
$c(\alpha)\restriction x;$ and vice versa,
 for every $\alpha\in\omega_1$ there is
$\beta\in\omega_1$ with $c(\beta)\restriction x$ eventually dominating
$b(\alpha)\restriction x.$
\enddefinition

It is immediate that locked sequences are $\leq_\bb$-equivalent.
Note that any bound on an infinite set $x\subset\omega$ of a collection
of increasing functions in $\rr$ easily yields a bound of that collection
on the whole $\omega.$

Now it is possible to lock unbounded sequences using one of the
standard tree forcings of \cite {BJ}:

\definition
{Definition 5.5} The Miller forcing $\Mi$ is the set of all nonempty
trees $T\subset{^{<\omega}\omega}$ consisting of increasing sequences
for which

\roster
\item for every $t\in T$ there is a splitnode $s$ of $T$ which extends $t$
\item if a sequence $s$ is a splitnode of $T$ then $s$ has in fact infinitely
many immediate successors in $T.$
\endroster

$\Mi$ is ordered by inclusion.
\enddefinition

The Miller forcing is proper, $\Cal M$-friendly--see Definition 6.6--and 
as such preserves
nonmeager sets of reals and unbounded sequences of functions
by the argument following Corollary 5.2. If $G\subset\Mi$
is a generic filter then $f=\bigcup\bigcap G\in\rr$ is an increasing function
called a Miller real.

\proclaim
{Theorem 5.6} $\Mi\Vdash$``every two unbounded sequences from the ground
model are locked''.
\endproclaim

\proclaim {Corollary 5.7}
It is consistent with ZFC that there are exactly two classes of 
$\leq_\bb$-equivalence.
\endproclaim

\demo {Proof}
Start with a model of ZFC+GCH and iterate Miller forcing $\omega_2$ times
with countable support. The resulting poset has $\aleph_2$-c.c.,
it is proper and $\Cal M$-friendly \cite {BJ}, therefore it does not
collapse unbounded sequences and forces $\bb=\aleph_1.$ In the resulting
model, there are exactly two classes of $\leq_\bb$-equivalence:
the objects which are not unbounded sequences and the unbounded sequences,
which are pairwise locked by the above theorem and a chain condition argument.
\qed
\enddemo

It also follows from the Theorem that whenever there are two unbounded 
sequences, one of length $\omega_1$ and the other of length $\omega_2,$
Miller forcing necessarily collapses $\aleph_2$ to $\aleph_1.$

\demo {Proof of Theorem 5.6}
Let $\dot x$ be an $\Mi$-name for the range of the Miller real. We shall show
that every two unbounded sequences $b,c:\omega_1\to\rr$ in the ground
model are forced to be locked by $\dot x.$ To this end, given
$T\in\Mi$ and $\alpha\in\omega_1$ a tree $S\in\Mi,S\subset T$ and an
ordinal $\beta\in\omega_1$ will be produced such that $S\Vdash$``$c(
\alpha)\restriction \dot x$ is eventually dominated by $b(\beta)
\restriction\dot x$''. The theorem then follows by the obvious density
and symmetricity arguments.

So fix $b,c,T$ and $\alpha$ as above. Let $\beta\in\omega$ be an ordinal 
such that $b(\beta)$ is not eventually dominated by any function recursive
in $c(\alpha)$ and $T.$ A tree $S\in\Mi,S\subset T$ with the same trunk
$t$ as $T$ will be found such that

$$
s\in S,n\in\dom(s)\setminus\dom(t)\text{ implies }c(\alpha)(n)\leq b(\beta)
(n). \tag {*}
$$

This will complete the proof. Let $S$ be defined by $s\in S$ iff $s\in T$ and
if $s'$ is the least splitnode of $T$ above or equal to $s$ then
for every $n\in\dom(s')\setminus\dom(t)$ it is the case that
$c(\alpha)(n)\leq b(\beta)(n).$

Obviously $t\in S\subset T$ and $S$ has property (*), moreover
$S$ is closed under initial segment and 
if $s\in S$ then the least splitnode of $T$ above or equal to $s$
belongs to $S$ as well. We must show that $S\in
\Mi,$ and this will follow from the fact that if $s\in S$ is a splitnode
of the tree $T$ then $s$ has infinitely many immediate successors
in $S.$ And indeed, let $y\subset\omega$ be the infinite set
of all integers $n\in\omega$ with $s^\smallfrown\langle n\rangle\in T$
and let $g:y\to\omega$ be a function defined by $g(n)=c(\alpha)(s'(m-1)),$
where $s'$ is the least splitnode of $T$ above or equal to
$s^\smallfrown\langle n\rangle$ and $m=\lth(s').$
Then by the choice of the ordinal $\beta\in\omega_1$ the set
$z=\{ n\in y:g(n)\leq b(\beta)(n)\}\subset\omega$ is infinite and every
sequence $s^\smallfrown\langle n\rangle:n\in z$ belongs to the tree $S.$ \qed
\enddemo

\subhead
{5.1. A model for $\bb=\aleph_1$}
\endsubhead

The $P_{max}$ variant for $\bb=\aleph_1$ will be built using the 
following notion of a witness: $b:\omega_1\to\rr$ is a good unbounded
sequence if it is unbounded and for every sequence $\eta\in{^{<\omega}\omega}$
the set $\{\alpha\in\omega_1:b(\alpha)$ pointwise dominates $\eta$
on its domain$\}\subset\omega_1$ is stationary. Note that if $j$ is a 
full iteration of a model $M,$ $M\models$``$b$ is a good unbounded sequence''
and $j(b)$ is unbounded then in fact $j(b)$ is a good unbounded sequence.
Also, whenever $\delta$ is a Woodin cardinal of $M$ and $j_\Bbb Q$ is the
$\Bbb Q_{<\delta}$-term for the canonical ultrapower embedding of $M$
then $j_\Bbb Q b(\omega_1^M)$ is a name for an unbounded function.

\proclaim
{Optimal Iteration Lemma 5.8} Assume $\bb=\aleph_1.$ Whenever $M$
is a countable transitive model of ZFC iterable with respect to
its Woodin cardinal $\delta$ and 
$M\models$``$b$ is a good unbounded sequence''
there is a full iteration $j$ of $M$ based on $\delta$ 
such that $j(b)$ is an unbounded
sequence.
\endproclaim

\demo {Proof}
Drawing on the asumption, choose an unbounded sequence $c$ of length $\omega_1$
and fix an arbitrary iterable model $M$ with $M\models$``$b$ is
an unbounded sequence and $\delta$ is a Woodin cardinal''. Two
full iterations $j_0,j_1$ of the model $M$ will be constructed simultaneously so
that the function $n\mapsto\max\{j_0b(\theta_{0,\alpha})(n),
j_1b(\theta_
{1,\alpha})(n)\}$ eventually dominates the function $c(\alpha),$ this for
every $\alpha\in\omega_1.$ Here $\theta_{0,\alpha}$ is $\omega_1$ in the
sense of the $\alpha$-th model on the iteration $j_0;$ similarly for
$\theta_{1,\alpha}.$

It follows immediately that one of the sequences $j_0(b),j_1(b)$ must
be unbounded, since if both were bounded--say by functions $f_0,f_1$
respectively--then the sequence $c$ would be bounded as well by the function
$n\mapsto\max\{f_0(n),f_1(n)\},$ contrary to the choice of $c.$

Now the iterations $j_0,j_1$ can be constructed easily using standard
bookkeeping arguments and the following claim. 

\proclaim {Claim 5.9}
Let $M_0,M_1$ be countable transitive models of ZFC and let

\roster
\item $M_0\models P$ is a poset, $p\in P,$ and $p\Vdash_P\dot x\in\rr$
is an increasing function unbounded over $M_0.$
\item $M_1\models Q$ is a poset, $q\in Q,$ and $q\Vdash_Q\dot y\in\rr$
is and increasing function unbounded over $M_1.$
\endroster

Suppose $f\in\rr$ is an arbitrary function. Then there are $M_0$ ($M_1,$
respectively) generic filters $p\in G\subset P,q\in H\subset Q$ 
such that the function
$n\mapsto\max\{(\dot x/G)(n),(\dot y/H)(n)\}$ eventually dominates $f.$
\endproclaim

\demo {Proof}
Without loss of generality assume that $f$ is an increasing function. 
Let $C_k:k\in\omega$ and $D_k:k\in\omega$ enumerate all open dense subsets of
$P$ in the model $M_0$ and of $Q$ in $M_1$ respectively. By induction
on $k\in\omega$ simultaneously build sequences $p=p_0\geq p_1\geq\dots\geq
p_k\geq\dots$ of conditions in $P,$ $q=q_0\geq q_1\geq\dots\geq q_k\geq\dots$
of conditions in $Q$ and integers $0=n_0=m_0, m_k\leq n_k<m_{k+1}$ so that

\roster
\item $p_{k+1}\in C_k,q_{k+1}\in D_k$ for all $k\in\omega$
\item $p_{k+1}$ decides $\dot x\restriction m_k,$ $q_{k+1}$ decides
$\dot y\restriction n_k$
\item for $k>0,$ for infinitely many integers $i\in\omega$ there is a condition
$r_i\leq p_k$ such that $r_i\Vdash_P$``$\dot x(m_k)=i$''; similarly,
 for infinitely many integers $i\in\omega$ there is a condition
$s_i\leq q_k$ such that $s_i\Vdash_Q\dot y(n_k)=i$
\item for $k>0,$ $p_{k+1}\Vdash\dot x(m_k)>f(n_k)$ and $q_{k+1}\Vdash\dot y
(n_k)>f(m_{k+1}).$
\endroster

This is easily arranged--(3) is made possible by the fact that $\dot x,\dot y$
are terms for unbounded functions. Let $G\subset P,H\subset Q$ be the filters
generated by the $p_i$'s or $q_i$'s respectively. These filters have
the desired genericity properties by (1) above and the function
$n\mapsto\max\{(\dot x/G)(n),(\dot y/H)(n)\}$ pointwise dominates $f$
from $n_1$ on. This can be argued from the induction hypothesis (4) and
the fact that $\dot x/G,\dot y/H$ and $f$ are all increasing functions. \qed
\enddemo

\enddemo
 
\proclaim 
{Strategic Iteration Lemma 5.10} Assume the Continuum Hypothesis.
The good player has a winning strategy
in the game $\Cal G_{\bb=\aleph_1}.$
\endproclaim

\demo {Proof}
Let $\vec N=\langle b,N_i,\delta_i:i\in\omega\rangle$ 
be a sequence of models with
a good unbounded sequence, let $y_0\in\Bbb Q_{\vec N}$ and let
$f\in\rr$ be an arbitrary function. We shall show that there is
an $\vec N$-generic filter $G\subset\Bbb Q_{\vec N}$ with $y_0\in G$
such that the function $j_{\Bbb Q}b(\omega_1^{\vec N})$ is
not eventually dominated by the function $f,$ where $j_\Bbb Q$ is
the $\Bbb Q_{<\delta_0}^{N_0}$-generic ultrapower of the model
$N_0$ using the filter $G.$ Then the winning strategy of the good
player consists essentially only of a suitable bookkeping using
the Continuum Hypothesis.

Now the proof of the existence of such a filter $G$ is in fact very easy.
We indicate a slightly inefficient though conceptual proof. Just
use the subgenericity Corollary 5.2 and the proof of Lemma 2.8 
with the following
changes:

\roster
\item $d$ is replaced with $b,$ the Hechler forcing is replaced with
Cohen forcing and the real $e\in\rr$ will be taken sufficiently
$\Bbb C$-generic
\item step (5) in the proof of Lemma 2.8 is replaced by: $e$ is not eventually
dominated by the function $f$
\item the poset $R_i$ is now the set of all pairs $\langle y,\eta\rangle$
with $y\in\Bbb Q_i,\eta\in{^{<\omega}\omega}$ and for every
$x\in y$ the sequence $\eta$ is everywhere on its domain dominated
by the function $b(x\cap\omega_1).$
\endroster

\qed
\enddemo

\proclaim {Conclusion 5.11}
The sentence $\bb=\aleph_1$ is $\Pi_2$-compact, even with the predicate 
$\frak B$ for unbounded sequences added to the language of $\stru.$
\endproclaim

\demo {Proof}
$\Pi_2$-compactness follows from Lemmas 5.8, 5.10.
 To show that $\Pi_2$ statements
of $\langle H_{\aleph_2},\in,\frak I,\frak B\rangle$ reflect to the model
$L(\Bbb R)^{P_{\bb=\aleph_1}},$ proceed as in Corollary 1.17.

Suppose in $V,$ a $\Pi_2$-sentence $\psi$--equal to $\forall x\ \exists y\
\chi(x,y)$ for some $\Sigma_0$ formula $\chi$--holds in
$\langle H_{\aleph_2},\in,\frak I,\frak B\rangle$ together with
$\bb=\aleph_1$ and let $\delta<\kappa$ be a Woodin and a measurable
cardinal respectively. For contradiction, suppose that $p\in P_{\bb=\aleph_1}$
forces $\lnot\psi$ to hold; strengthening the condition $p$ if necessary
we may assume that for some $x\in M_p,$ $p\Vdash\forall y\lnot\chi
(\dot k_p(x),y)$ where $\dot k_p$ is the term for the canonical iteration
of $M_p$ as defined in 1.14. By Corollary 1.8, there is a countable
transitive  model $M$ elementarily embeddable into $V_\kappa$
such that $p\in M$ and $M$ and all of its generic extensions by
posets of size $\leq\frak c^M$ are iterable. So $M\models$``$\bb=\aleph_1$
and $\langle H_{\aleph_2},\in,\frak I,\frak B\rangle\models\psi$ and some
cardinal $\bar\delta$ is Woodin''. The optimal iteration lemma
applied in $M$ yields there a full iteration $j$ of the model $M_p$
such that $j(b_p)$ is a good unbounded sequence. Let $N$ be an
$\Mi$-generic extension of the model $M_p$ and let
$q=\langle N,j(b_p),\bar\delta,j(H_p)\cup\{\langle j,p\rangle\}\rangle.$

First, $q\in P_{\bb=\aleph_1}$ is a condition strengthening $p.$ To see that
note that the model $N$ is iterable, the iteration $j$ is full in $N,$
since the forcing $\Mi$ preserves stationary sets, and $j(b_p)$ is
a good unbounded sequence in $N$ since $\Mi$ preserves such sequences.

Second, $q\Vdash$``$\frak B$ in the sense of the structure 
$\dot k_q(H_{\aleph_2})^M$ is just $\frak B\cap\dot k_q(H_{\aleph_2})^M$''.
Observe that $N\models$``all sequences in $(H_{\aleph_2})^M$ which
there are unbounded are locked with $j(b_p)$'', by Theorem 5.6
applied in the model
$M.$ Therefore $q$ forces even the following stronger statement,
by the elementarity of the embedding $\dot k_q:$
``whenever $\dot k_q(H_{\aleph_2})^M\models$`$c$ is an unbounded sequence'
then $c$ and $\dot k_qj(b_p)$ are locked; since $\dot k_qj(b_p)$
is unbounded, the sequence $c$ must be unbounded as well''.

Third, $q\Vdash$``$\exists y\ \chi(\dot k_qj(x)=\dot k_p(x),y)$'', giving
the final contradiction with our choice of $p$ and $x.$ Since $\psi$
holds in the model $M,$ there must be some $y\in(H_{\aleph_2})^M$ such that
$\langle H_{\aleph_2},\in,\frak I,\frak B\rangle^M\models\chi(j(x),y).$
Then $q\Vdash\langle H_{\aleph_2},\in,\frak I,\frak B\rangle^M\models
\chi(\dot k_qj(x),\dot k_q(y))$ by the previous paragraph and absoluteness
of $\Sigma_0$ formulas. \qed
\enddemo

\head
{6. Uniformity of the meager ideal}
\endhead

The sentence ``there is a nonmeager set of reals of size $\aleph_1$''
does not seem to be compact, however, a similar a bit stronger
assertion is.

\definition
{Definition 6.1} A sequence $\langle r_\alpha:\alpha\in\omega_1\rangle$ of real
numbers is called weakly Lusin if for every meager $X\subset\Bbb R$
the set $\{\alpha\in\omega_1:r_\alpha\in X\}\subset\omega_1$ is nonstationary.
\enddefinition

Thus the existence of a weakly Lusin sequence is a statement intermediate
between a nonmeager set of size $\aleph_1$ and a Lusin set. It is equivalent
to neither of them, as will be shown below.
\subhead
{6.0. A model for a weakly Lusin sequence}
\endsubhead

We will prove the iteration lemmas necessary to conclude that the sentence
$\phi=$``there is a weak Lusin sequence'' is $\Pi_2$-compact. Note that if
$\langle M,k,\delta\rangle$ is a triple such that $M\models$``$k$
is a weak Lusin sequence and $\delta$ is a Woodin cardinal'' then in $M,$
$\Bbb Q_{<\delta}\Vdash$``$j_{\Bbb Q}k(\omega_1^M)$ 
is a Cohen real over $M",$
and so it generates a natural Cohen subalgebra of $\Bbb Q_{<\delta}.$ The
following abstract copying lemma will be relevant:

\proclaim {Lemma 6.2}
Let $N$ be a countable transitive model of ZFC, $N\models$``$P$ is
a partially ordered set and $P\Vdash\dot r$ is a Cohen real''. Suppose
that $p\in P$ and $s\in\rr$ is a Cohen real over $N.$ Then there is
an $N$-generic filter $G\subset P$ so that $p\in G$ and $\dot r/G
=s$ modulo finite.
\endproclaim

\demo {Proof}
In the model $N,$ let $P=\Bbb C*\dot Q$ where $\Bbb C$ is the cohen
subalgebra of $\Bbb B$ generated by the term $\dot r.$ It follows from
the assumptions and some boolean algebra theory in $N$ that there are
$q\leq p$ in $RO(P)$ and a finite sequence $\eta$ such that $pr_{\Bbb C}(q)=
[[\eta\subset\dot r]]_{\Bbb C}.$ Let $t\in\rr$ be the function defined by
$\eta\subset t$ and $s(n)=t(n)$ for $n\neq dom(\eta).$ Since the real
$s$ is Cohen over $N$ and $s=t$ modulo finite, even $t$ is Cohen and
the filter $H\subset\Bbb C$ generated by the equation $\dot r=t$
is $N$-generic. Finally, choose an $N$-generic filter $G\subset P$
with $H\subset G$ and $q\in G.$ This is possible since $pr_{\Bbb C}(q)\in H.$
The filter $G\subset P$ is as desired. \qed
\enddemo

\proclaim {Optimal Iteration Lemma 6.3}
 Assume that there is a weakly Lusin sequence. Whenever
$M$ is a countable transitive model iterable with respect to its
Woodin cardinal $\delta$ such that $M\models$
``$k$ is a weakly Lusin sequence'',
then there is a full iteration $j$ of $M$ such that $j(k)$
is a weakly Lusin sequence.
\endproclaim

\demo {Proof}
Fix a Lusin sequence $J$ and $M,k,\delta$ as in the Lemma. 
We shall produce a full
iteration $j$ of $M$ such that there is a club $C\subset\omega_1$
with $\forall\alpha\in C\ jk(\alpha)=J(\alpha)$ modulo finite. Then
$jk$ really is a Lusin sequence and the Lemma is proved.

The desired iteration $j$ will be constructed by induction on $\alpha\in
\omega_1.$ Let $S_\xi:\xi\in\omega_1$ be a partition of $\omega_1$ into 
pairwise disjoint stationary sets. By induction on $\alpha\in\omega_1,$
 models $M_\alpha$ together with the elementary embeddings will be built.
plus an enumeration $\{\langle x_\xi,\beta_\xi\rangle:\xi\in\omega_1\}$
of all pairs $\langle x,\beta\rangle$ with $x\in \Bbb Q_\beta.$
The induction hypotheses at $\alpha$ are:

\roster
\item if $\gamma<\alpha$ and $J(\theta_\gamma)$ is a Cohen real
over $M_\gamma$ then $k_{\gamma+1}(\theta_\gamma)=J(\theta_\gamma)$
modulo finite
\item if $\gamma\leq\alpha$ then $\{\langle x_\xi,\beta_\xi\rangle:
\xi\in\theta_\gamma\}$ enumerates all pairs
 $\langle x,\beta\rangle$ with $\beta<\gamma,$
 $x\in \Bbb Q_\beta$
\item if $\gamma<\alpha$ and $\theta_\gamma\in S_\xi$ for some $\xi\in
\theta_\gamma$ then $j_{\beta_\xi,\gamma}(x_\xi)\in G_\gamma.$
\endroster

As before, the hypothesis (1) makes sure that the sequence $J$ gets copied
onto $jk$ properly and (2,3) are just bookkeping tools for making
the resulting iteration full.

At limit ordinals just direct limits are taken and the new enumeration is
the union of all old ones. The successor step is handled easily using
the previous Lemma applied for $N=M_\alpha,P=\Bbb Q_\alpha,
\dot r=j_{\Bbb Q}k_\alpha(\theta_\alpha),s=J(\theta_\alpha)$ and $
b=j_{\beta_\xi\alpha}(x_\xi),$ the last two in the case that
$s$ is Cohen over $M_\alpha$ and $\theta_\alpha\in S_\xi$ for
some (unique) $\xi\in\theta_\alpha.$

To prove that the resulting iteration is as desired, note that it is
full and that the set $D=\{\alpha\in\omega_1:J(\theta_\alpha)$
is a Cohen real over $M_\alpha\}$ contains a closed unbounded set.
For assume otherwise. Then the complement $S$ of $D$ is stationary
and for every limit ordinal $\alpha\in S$ there is a nowhere dense
tree in some $M_\beta,\beta\in\alpha$ such that the real
$J(\omega_1^{M_\alpha})$ is a branch of this tree--this is because
a direct limit is taken at step $\alpha.$ By a simple Fodor-style
argument there is a nowhere dense tree and a stationary set $T\subset S$
such that every $J(\theta_\alpha):\alpha\in T$ is a branch of this
tree. This contradicts the assumption of $J$ being a weakly Lusin sequence.
\qed
\enddemo

\proclaim {Strategic Iteration Lemma 6.4} Assume the Continuum Hypothesis.
The good player has a winning strategy
in the game $\Cal G_{\phi}$ connected with weakly Lusin sequences.
\endproclaim

\demo {Proof}
Given a sequence $\vec N=\langle k,N_i,\delta_i:i\in\omega\rangle$
of models with a weakly Lusin sequence $k,$ a condition
$y_0\in\Bbb Q_{\vec N}$ and nowhere dense trees $T_n:n\in\omega,$
we shall show that there is an $\vec N$-generic filter $G\subset
\Bbb Q_{\vec N}$ with $y_0\in G$ such that the real
$j_{\Bbb Q}k(\omega_1^{\vec N})$ is not a branch through any of the
trees $T_n,$ where $j_\Bbb Q$ is the generic ultrapower embedding
of the model $N_0$ using the filter $G\cap\Bbb Q_{<\delta_0}^{N_0}.$
With this fact in hand, a winning strategy for the good player
consists of just a suitable bookkeping using the Continuum Hypothesis.

Again, we provide maybe a little too conceptual proof of the existence
of the filter $G,$ using the ideas from Lemma 2.8. No subgenericity 
theorems are needed this time. Let $i\in\omega$ and work in $N_i.$
Let $\Bbb Q_i=\Bbb Q_{\delta_i}$ and $j_i$ to be the $\Bbb Q_i$
term for the generic ultrapower embedding of the model $N_i.$
So the term $j_ik(\omega_1^{\vec N})$ is a term for a Cohen real,
and it gives a complete embedding of $\Bbb C$ into $\Bbb Q_i.$
The key point is that with this embedding, the computation of
the projection $pr_{\Bbb C}(y)$ gives the same value in $\Bbb C$
in every model $N_i$ with $y\in \Bbb Q_i,$ namely
$\Sigma_{\Bbb C}\{\eta\in{^{<\omega}\omega}:$ the system $\{x\in y:
\eta\subset k(x\cap\omega_1)\}$ is stationary$\}.$

First, fix a suitably generic filter $H\subset \Bbb C.$ The requirements
are:

\roster
\item $pr_\Bbb C(y_0)\in H$
\item $H$ meets all the maximal antichains that happen to belong to
$\bigcup_i N_i$
\item the Cohen real $c\in \rr$ given by the filter $H$ does not
constitute a branch through any of the nowhere dense trees $T_n.$
\endroster

This is easily done. Now let $X_n:n\in\omega$ be an enumeration
of all maximal antichains of $\Bbb Q_{\vec N}$ in $\bigcup N_i$
and by induction on $n\in\omega$ build a decreasing
sequence $y_n:n\in\omega$ of conditions in $\Bbb Q_{\vec N}$
so that

\roster
\item $pr_\Bbb C(y_n)\in H$
\item $y_{n+1}$ has an element of $X_n$ above it.
\endroster

This can be done since the filter $H$ is $\Bbb C$-generic over every 
model $N_i.$ In the end, let $G$ be the filter on $\Bbb Q_{\vec N}$
generated by the conditions $y_n:n\in\omega.$ This filter is as desired.
Note that $c$ is the uniform value of $j_ik(\omega_1^{\vec N})$
as evaluated according to this filter. \qed
\enddemo

\proclaim {Conclusion 6.5}
The sentence ``there is a weakly Lusin sequence of reals'' is $\Pi_2$-compact.
\endproclaim

It is unclear whether it is possible to add a predicate for witnesses
in this case.
\subhead
{6.1. Combinatorics of weakly Lusin sequences}
\endsubhead

In this subsection it is proved that the existence of a Lusin set,
weakly Lusin sequence and a nonmeager set of size $\aleph_1$ are nonequivalent
assertions.

The following regularity property of forcings will be handy:

\definition
{Definition 6.6} 
\cite {BJ, 6.3.15} A forcing $P$ is called $\Cal M$-friendly if for every
large enough regular cardinal $\lambda,$ every condition $p\in P,$
every countable elementary submodel $M$ of $H_\lambda$ with $p,P\in M$ 
and every function $h\in\rr$ Cohen-generic over $M$ there is a
condition $q\leq p$ such that $q$ is master for $M$ and $q\Vdash$
``$\check h$ is Cohen-generic over $M[G]$''.
\enddefinition

It is not difficult to prove that $\Cal M$-friendly forcings preserve
nonmeager sets. Moreover $\Cal M$-friendliness is preserved under
countable support iterations \cite {BJ, Section 6.3.C}. 
The following fact, pointed out to us
by Tomek Bartoszy\' nski, replaces our original more complicated argument.

\proclaim
{Lemma 6.7}\roster
\item The Miller forcing $\Mi$--see Definition 5.5--is $\Cal M$-friendly.
\item $\Mi$ destroys all weakly Lusin sequences from the ground model.
\endroster
\endproclaim

\proclaim {Corollary 6.8}
It is consistent with ZFC that there is a nonmeager set of size $\aleph_1$
but no weakly Lusin sequences.
\endproclaim

\demo {Proof}
Iterate Miller forcing over a model of GCH $\omega_2$ times. \qed
\enddemo

\demo {Proof of Lemma}
Suppose $M$ is a countable transitive model of a rich fragment of ZFC,
$T\in\Mi\cap M$ and $f\in\rr$ is a Cohen real over $M.$ We shall produce
$M$-master conditions $S_0,S_1\subset T$ in $\Mi$ such that

\roster
\item $S_0\Vdash\check f$ is a Cohen real over $M$
\item $S_1\Vdash\check f$ is eventually dominated by the Miller real.
\endroster

This will finish the proof: (1) shows that $\Mi$ is $\Cal M$-friendly
and (2) by a standard argument implies that for any weakly Lusin
sequence $\langle r_\alpha:\alpha\in\omega_1\rangle$ of elements of
$\rr$ the set $\{\alpha\in\omega_1:r_\alpha$ belongs to the meager
set of all functions in $\rr$ eventually dominated by the Miller real
$\}\subset\omega_1$ is $\Mi$-forced to be stationary.

By a mutual genericity argument there is an $M$-generic filter
$G\subset Coll(\omega,(2^\frak c)^M)$ such that the function
$f$ is still Cohen generic over $M[G].$ Work in $M[G].$ Let $\eta_k:k\in\omega$
enumerate the Cohen forcing ${^{<\omega}\omega},$ let $\dot O_k:
k\in\omega$ enumerate all the $\Mi\cap M$-names for dense subsets of
the Cohen forcing in $M$ and let $D_k:k\in\omega$ enumerate the open
dense subsets of $\Mi\cap M$ in $M.$ Build a fusion sequence
$T=T_0\geq T_1\geq T_2\geq\dots$ of trees in $\Mi\cap M$ so that if $s$ is a $k$-th level splitnode of $T_k$ and $\{i_n:n\in\omega\}$
is an enumeration of the set of all integers $i$ with $s^\smallfrown
\langle i\rangle\in T_k$ then for every $n\in\omega$ the sequence
 $s^\smallfrown\langle i_n\rangle$ belongs to $T_{k+1},$ $T_{k+1}\restriction
 s^\smallfrown\langle i_n\rangle\in D_k$ and there is some extension
$\eta$ of $\eta_n$ in the Cohen forcing ${^{<\omega}\omega}$ such that
$T_{k+1}\Vdash\check\eta\in\dot O_k.$
This is readily done. Let $S=\bigcup_i T_i\in\Mi\cap M[G].$ The tree $S$
is an $M$-master condition in $\Mi$ and since $f\in\rr$ is Cohen generic
over the model $M[G],$ the following two subtrees $S_0,S_1$ of $S$
are still in $\Mi:$

\roster
\item $S_0=\{s\in S:$ whenever $s$ is a proper extension of a $k$-th
level splitnode $t\in T_k$ then $T_{k+1}\restriction s\Vdash\check\eta\in
\dot O_k$ for some $\eta\subset f\}$
\item $S_1=\{s\in S:\forall n\in dom(s\setminus$the trunk of $T)\ f(n)<s(n)\}.$
\endroster

It is easy to see that the trees $S_0,S_1$ are as desired. \qed
\enddemo

\proclaim
{Question 6.9} Does the saturation of the nonstationary ideal plus
the existence of a nonmeager set of reals of size $\aleph_1$ imply
the existence of a weakly Lusin sequence?
\endproclaim

Next it will be proved that the existence of a weakly Lusin sequence does
not imply that of a Lusin set. A classical forcing argument can be tailored
to fit this need; instead, we shall show that there are no Lusin
sets in the model built in the previous Subsection. This follows from
Theorem 1.15(2) and a simple density argument using the following fact:

\proclaim {Lemma 6.10}
($AD^{L(\Bbb R)}$) Let $M$ be a countable transitive iterable model,
$M\models$``$K$ is a weakly Lusin sequence and $X$ is a Lusin set''.
Then there is a countable transitive iterable model $N$ and an iteration
$j\in N$ such that
$N\models$``$j$ is a full iteration of $M,$ $j(K)$ is a weakly Lusin
sequence and $j(X)$ is not a Lusin set''.
\endproclaim

\demo {Proof}
Fix $M,K$ and $X$ and for notational reasons assume that $X$ is really
an injective function from $\omega_1^M$ to $\rr$ enumerating that Lusin set.
Working in the model $M,$ if $\delta$ is a Woodin cardinal, $\Bbb Q_{<\delta}$
is the nonstationary tower and $j_{\Bbb Q}$ is the $\Bbb Q_{<\delta}$-name
for the canonical embedding of $M$ then both $j_{\Bbb Q}K(\omega_1^M)$
and $j_{\Bbb Q}X(\omega_1^M)$ are $\Bbb Q_{<\delta}$ terms for Cohen reals.

Using the determinacy assumption, choose a countable transitive model
$N_0$ such that $M\in N_0$ is countable there and $N_0$ and all of its
generic extensions by forcings of size $\aleph_1^{N_0}$ are iterable.
Work in $N_0.$ Force two $\omega_1$-sequences $\langle c_\beta:\beta
\in\omega_1\rangle,\langle d_\beta:\beta\in\omega_1\rangle$ of
Cohen reals--elements of $\rr$--and a function $e\in\rr$ eventually dominating
every $d_\beta:\beta\in\omega_1$ with finite conditions.

Set $N=N_0[\langle c_\beta:\beta
\in\omega_1\rangle,\langle d_\beta:\beta\in\omega_1\rangle,e]$ In the model
$N,$ the reals $\langle c_\beta:\beta\in\omega_1\rangle$ constitute
a weakly Lusin sequence, indeed a Lusin set, because their sequence is
Cohen generic over the model $N_0[\langle d_\beta:\beta\in\omega_1\rangle,
e].$ In the model $N_0[\langle c_\beta:\beta
\in\omega_1\rangle,\langle d_\beta:\beta\in\omega_1\rangle]$ build a full
iteration $j$ of $M$ so that

\roster
\item if $\alpha\in\omega_1$ is limit then whenever made possible by
the model $M_\alpha$ we have $jK(\theta_\alpha)=c_{\theta_\alpha}$
modulo finite
\item if $\alpha\in\omega_1$ is successor then $jX(\theta_\alpha)$
is equal to one of the reals $d_\beta:\beta\in\omega_1$ modulo finite.
\endroster

By the arguments from the previous subsection and the fact that
$\langle c_\beta:\beta
\in\omega_1\rangle,\langle d_\beta:\beta\in\omega_1\rangle$ are
mutually generic sequences of Cohen reals over $N_0,$ (2) is always possible
to fulfill and the set $\{\beta\in\omega_1:jK(\beta)=c_\beta$ modulo finite
$\}\subset\omega_1$ will contain a club.

Now $N,j$ are as desired. In the model $N,$ the iteration $j$ is full since
$N$ is a c.c.c. extension of $N_0[\langle c_\beta:\beta
\in\omega_1\rangle,\langle d_\beta:\beta\in\omega_1\rangle]$ in which
$j$ was constructed to be full; the sequence $jK$ is on a club modulo finite
equal to a weakly Lusin sequence $\langle c_\beta:\beta
\in\omega_1\rangle$ and so is weakly Lusin iself; and the set $\{jX
(\theta_\alpha):\alpha\in\omega_1$ successor$\}$ is an uncountable subset
of $\rng(jX)$ contained in the meager set of all reals eventually dominated by
$e\in\rr,$ consequently $\rng(jX)$ is not a Lusin set. \qed
\enddemo
 
\subhead
{6.2. The null ideal}
\endsubhead

The methods of this paper can be adapted to give a parallel result about 
the null ideal.

\definition {Definition 6.11}
A sequence $\langle r_\alpha:\alpha\in\omega_1\rangle$ of real numbers
is called a weakly Sierpi\' nski sequence if for every null set $S,$
the set $\{\alpha\in\omega:r_\alpha\in S\}$ is nonstationary.
\enddefinition

\proclaim {Theorem 6.12}
The sentence ``there is a weakly Sierpi\' nski sequence'' is
$\Pi_2$-compact.
\endproclaim

By an argument parallel to 6.7(2) it can be proved that existence
of a nonnul set of size $\aleph_1$ and of a weakly Sierpi\' nski sequence
are nonequivalent statements.

%\input refs
%refs


\Refs\widestnumber\key{BJZ}
\ref
 \key AS
 \by U. Abraham and S. Shelah
 \paper A $\Delta_2^2$ wellorder of the reals and incompactness
of $L(Q^{MM})$
 \jour Ann. Pure Appl. Logic
 \vol 59
 \yr 1993
 \pages 1--32
\endref
\ref
 \key Ba
 \by T. Bartoszy\' nski
 \paper Combinatorial aspects of measure and category
 \jour Fund. Math.
 \vol 127
 \yr 1987
 \pages 225--239
\endref
\ref
 \key BJ 
 \by T. Bartoszy\' nski and H. Judah
 \book Set theory. On the structure of the real line
 \publ A K Peters
 \publaddr Wellesley, MA
 \yr 1995
\endref
\ref
 \key BJZ
 \by B. Balcar, T.J. Jech and J. Zapletal
 \paper Generalizations of Cohen algebras
 \paperinfo preprint
 \yr 1995
\endref
\ref
 \key J1
 \by T. J. Jech
 \paper Nonprovability of Souslin's Hypothesis
 \jour Comment. Math. Univ. Carolinae
 \vol 8
 \yr 1967
 \pages 291--305
\endref
\ref
 \key J2
 \bysame
 \book Set Theory
 \publ Academic Press
 \publaddr New York
 \yr 1978
\endref
\ref
 \key JS
 \by H. Judah and S. Shelah
 \paper Suslin forcing
 \jour J. Symbolic Logic
 \vol 53
 \yr 1988
 \pages 1188--1207
\endref
\ref
 \key K
 \by S. Koppelberg
 \paper Characterization of Cohen algebras
 \inbook Paper on General Topology and Applications
 \eds S. Andima, R. Kopperman, P. R. Misra and A. R. Todd
 \bookinfo Annals of the New York Academy of Sciences 704
 \yr 1993
 \pages 227--237
\endref
\ref
 \key L
 \by R. Laver
 \paper On the consistency of Borel's conjecture
 \jour Acta Math.
 \vol 137
 \yr 1976
 \pages 151--169
\endref
\ref
 \key S1
 \by S. Shelah
 \book Proper forcing
 \publ Springer Verlag
 \bookinfo Lecture Notes in Math. 940
 \publaddr New York
 \yr 1982
\endref
\ref
 \key S2
 \bysame
 \paper Can you take Solovay inaccessible away?
 \jour Israel J. Math.
 \vol 48
 \yr 1984
 \pages 1--47
\endref
\ref
 \key S3
 \bysame
 \paper Vive la difference I
 \inbook Set Theory of the continuum
 \eds H. Judah, W.Just and W. H. Woodin
 \yr 1992
 \pages 357--405
\endref
\ref
 \key SZ
 \by S. Shelah and J. Zapletal
 \paper Embeddings of Cohen algebras
 \jour Adv. Math.
 \toappear
\endref
\ref
 \key So
 \by M. Souslin
 \paper Probl\` eme 3
 \jour Fund. Math.
 \vol 1
 \yr 1920
 \pages 223
\endref
\ref
 \key Te
 \by S. Tennenbaum
 \paper Souslin's problem
 \jour Proc. Nat. Acad. Sci. USA
 \vol 59
 \yr 1968
 \pages 60--63
\endref
\ref
 \key T1
 \by S. Todorcevic
 \paper Partitioning pairs of countable ordinals
 \jour Acta Math.
 \vol 159
 \yr 1987
 \pages 261--294
\endref
\ref
 \key T2
 \bysame
 \book Partition problems in topology
 \bookinfo Cont. Math. 84
 \publ Amer. Math. Soc.
 \publaddr Providence
 \yr 1989
\endref
 
\ref
 \key Tr
 \by J. K. Truss
 \paper Connections between different amoeba algebras
 \jour Fund. Math.
 \vol 130
 \yr 1988
 \pages 137--155
\endref
\ref
 \key W1
 \by W. H. Woodin
 \paper Supercompact cardinals, sets of reals, and weakly homogeneous trees
 \jour Proc. Natl. Acad. Sci. USA
 \vol 88
 \yr 1988
 \pages 6587--6591
\endref
\ref
 \key W2
 \bysame
 \book The axiom of determinacy, forcing axioms, and the nonstationary ideal
 \toappear
\endref
\endRefs
\enddocument